\documentclass[12pt]{amsart}

\usepackage{amsmath}
\usepackage{amscd}
\usepackage{amssymb}
\usepackage{latexsym}
\usepackage{stmaryrd}
\usepackage{mathrsfs}

\usepackage[all]{xy}
\DeclareMathAlphabet{\mathbbb}{U}{bbold}{m}{n}

\newtheorem{theorem}{Theorem}[section]
\newtheorem*{theorem*}{Theorem}
\newtheorem{lemma}[theorem]{Lemma}
\newtheorem*{lemma*}{Lemma}
\newtheorem{corollary}[theorem]{Corollary}
\newtheorem{proposition}[theorem]{Proposition}

\newtheorem{remark}[theorem]{Remark}
\newtheorem{definition}[theorem]{Definition}

\numberwithin{equation}{section}

\makeatletter
\def\revddots{\mathinner{\mkern1mu\raise\p@
\vbox{\kern7\p@\hbox{.}}\mkern2mu
\raise4\p@\hbox{.}\mkern2mu\raise7\p@\hbox{.}\mkern1mu}}
\makeatother % ddots gespiegelt
\newcommand{\bgl}{\begin{equation}} %eine Gleichung mit Ziffer
\newcommand{\egl}{\end{equation}}
\newcommand{\bgloz}{\begin{equation*}} %eine Gleichung ohne Ziffer
\newcommand{\egloz}{\end{equation*}}
\newcommand{\bgln}{\begin{eqnarray}} %mehrere Gleichungen mit Ziffer
\newcommand{\egln}{\end{eqnarray}}
\newcommand{\bglnoz}{\begin{eqnarray*}} %mehrere Gleichungen ohne Ziffer
\newcommand{\eglnoz}{\end{eqnarray*}}
\newcommand{\btheo}{\begin{theorem}}
\newcommand{\etheo}{\end{theorem}}
\newcommand{\btheooz}{\begin{theorem*}}
\newcommand{\etheooz}{\end{theorem*}}
\newcommand{\blemma}{\begin{lemma}}
\newcommand{\elemma}{\end{lemma}}
\newcommand{\blemmaoz}{\begin{lemma*}}
\newcommand{\elemmaoz}{\end{lemma*}}
\newcommand{\bproof}{\begin{proof}}
\newcommand{\eproof}{\end{proof}}
\newcommand{\bbew}{\begin{beweis}}
\newcommand{\ebew}{\end{beweis}}
\newcommand{\bremark}{\begin{remark}\em}
\newcommand{\eremark}{\end{remark}}
\newcommand{\bdefin}{\begin{definition}}
\newcommand{\edefin}{\end{definition}}
\newcommand{\bprop}{\begin{proposition}}
\newcommand{\eprop}{\end{proposition}}
\newcommand{\bcor}{\begin{corollary}}
\newcommand{\ecor}{\end{corollary}}
\newcommand{\bfa}{\begin{cases}} %Fallunterscheidung
\newcommand{\efa}{\end{cases}}
\newcommand{\cD}{\mathcal D}
\newcommand{\cE}{\mathcal E}

\newcommand{\cG}{\mathcal G}
\newcommand{\cH}{\mathcal H}

\newcommand{\cK}{\mathcal K}
\newcommand{\cL}{\mathcal L}

\newcommand{\cO}{\mathcal O}

\newcommand{\cQ}{\mathcal Q}

\newcommand{\cT}{\mathcal T}
\newcommand{\cU}{\mathcal U}

\def\Az{\mathbb{A}}
\def\Cz{\mathbb{C}}

\def\Nz{\mathbb{N}}

\def\Qz{\mathbb{Q}}
\def\Rz{\mathbb{R}}
\def\Tz{\mathbb{T}}
\def\Zz{\mathbb{Z}}

\def\1z{\mathbbb{1}}
\newcommand{\fA}{\mathfrak A}

\def\xf{{\bf{x}}}
\def\yf{{\bf{y}}}
\def\zf{{\bf{z}}}
\newcommand{\an}[1]{``#1''} % Anfuehrungsstriche
\newcommand{\ti}{\tilde}

\newcommand{\ri}{\rightarrow}
\newcommand{\lori}{\longrightarrow}

\newcommand{\ma}{\mapsto} % wird abgebildet auf
\newcommand{\mafr}{\mapsfrom} % wird abgebildet auf
\newcommand\onto{\twoheadrightarrow} % surjektiv
\newcommand\into{\hookrightarrow} % injektiv
\def\SEMI{\mbox{$\times\kern-2pt\vrule height5pt width.6pt \kern3pt $}}

\newcommand{\halb}{\tfrac{1}{2}}

\newcommand{\End}{{\rm End}\,}

\newcommand{\Spec}{{\rm Spec\,}} % Spektrum

\newcommand{\id}{{\rm id}}

\newcommand{\staralg}{{\text{*-} \rm alg}}
\newcommand{\Ind}{\mathrm{ Ind}\,}

\newcommand{\Ad}{{\rm Ad\,}}

\newcommand{\img}{{\rm im\,}}

\newcommand{\reg}{^\times} % regulaer
\newcommand{\pos}{_{>0}} % positiv
\newcommand{\lspan}{{\rm span}} % linearer Aufspann
\newcommand{\abs}[1]{\lvert#1\rvert} % Betrag
\newcommand{\norm}[1]{\left\|#1\right\|} % Norm
\newcommand{\defeq}{\mathrel{:=}} % per Definition
\newcommand{\dop}{\text{: }} % in Mengen
\newcommand{\fuer}{\text{ for }} % bei Fallunterscheidungen
\newcommand{\fa}{\text{ for all }} % fuer alle
\newcommand{\ilim}{\varinjlim} % induktiver Limes
\newcommand{\plim}{\varprojlim} % projektiver Limes
\newcommand{\e}{{\rm e}} % exp(2 PI . )
\newcommand{\dotcup}{\ensuremath{\mathaccent\cdot\cup}} % disjunkte Vereinigung
\newcommand{\Ell}{{\rm L}} % L
\newcommand{\rte}{\overset{e}{\rtimes}} % verschraenktes Produkt mit Endomorphismen
\newcommand{\lge}{\left\{} % links geschweift
\newcommand{\rge}{\right\}} % rechts geschweift
\newcommand{\lru}{\left(} % links rund
\newcommand{\rru}{\right)} % rechts rund
\newcommand{\leck}{\left[} % links eckig
\newcommand{\reck}{\right]} % rechts eckig
\newcommand{\lsp}{\left\langle} % links spitz
\newcommand{\rsp}{\right\rangle} % links spitz
\newcommand{\rukl}[1]{\lru #1 \rru} % runde Klammer
\newcommand{\eckl}[1]{\leck #1 \reck} % eckige Klammer
\newcommand{\gekl}[1]{\lge #1 \rge} % geschweifte Klammer
\newcommand{\spkl}[1]{\lsp #1 \rsp} % spitze Klammer
\newcommand{\menge}[2]{\gekl{ #1 \dop #2 }} % Menge
\begin{document}

\title[The 2-adic ring $C^*$-algebra of the integers]{The 2-adic ring $C^*$-algebra of the integers and its representations}
\author[N. S. Larsen]{Nadia S. Larsen}
\address{Department of Mathematics, University of Oslo, PO BOX 1053 Blindern, N-0316 Oslo, Norway.}
\email{nadiasl@math.uio.no}
\author[X. Li]{Xin Li}
\address{Mathematisches Institut, Fachbereich Mathematik und Informatik der Universit{\"a}t M{\"u}nster, Einsteinstra{\ss}e 62, 48149 M{\"u}nster, Germany}
\email{xinli.math@uni-muenster.de}

\thanks{This research was supported by the Research Council of Norway and the Deutsche Forschungsgemeinschaft.}


\begin{abstract}
We study the 2-adic version of the ring $C^*$-algebra of the integers.
First, we work out the precise relation between the Cuntz algebra $\cO_2$ and our 2-adic ring $C^*$-algebra in terms of representations.
Secondly, we prove a 2-adic duality theorem identifying the crossed product arising from 2-adic affine transformations on the 2-adic numbers
with the analogous crossed product algebra over the real numbers.
And finally, as an outcome of this duality result,
we construct an explicit imprimitivity bimodule and prove that it transports
one canonical representation into the other.
\end{abstract}
\date{November 2, 2011}
\maketitle

\tableofcontents

\setlength{\parindent}{0pt} \setlength{\parskip}{0.5cm}

\section*{Introduction}

Recently Cuntz introduced a $C^*$-algebra $\cQ_\Nz$ attached to the $ax+b$-semigroup over the natural numbers in \cite{Cun2}. He proved,
among other things, that $\cQ_\Nz$ is
purely infinite, simple, and Morita equivalent to the crossed product of the algebra of continuous functions, vanishing at infinity, on the space of finite adeles
by the affine action of the $ax+b$-semigroup over $\Qz$. The algebra $\cQ_\Nz$ is universal
for a family of isometries $s_n$ for $n\in \Nz^\times$ and a single unitary $u$ which satisfy relations that encode the
multiplicative action of the non-zero natural numbers and the
additive  action of the natural numbers. The $ax+b$-semigroup of the natural numbers is a semidirect product $\Zz\rtimes\Nz^\times$
 coming from the action of $\Nz^\times$ on $\Zz$ given by multiplication, and replacing $\Nz^\times$ by $\Zz^\times$ gives rise to
 the $C^*$-algebra $\cQ_\Zz$ of the ring $\Zz$. A far-reaching generalisation of $\cQ_\Zz$ was obtained in \cite{CuLi1}, where the
 authors introduced a  $C^*$-algebra $\fA[R]$ of an integral domain
$R$, and unveiled interesting connections with algebraic number theory. In \cite{CuLi2},
the $K$-theory of algebras of the kind $\fA[R]$ was computed by means of a duality theorem.
The construction of $\fA[R]$ was generalised to a large class of rings by the second named author in \cite{Li}.

Our goal here is to show that the analogue of $\cQ_\Nz$ where the single prime $2$ acts by multiplication
on $\Zz$ is a $C^*$-algebra which shares many of the structural properties of  $\cQ_\Nz$, although it is
not a special case of $\fA[R]$ for any ring $R$ (at least not directly).
Moreover, we study the representations of the 2-adic analogue of $\cQ_\Nz$
and present its precise relation to the Cuntz algebra $\cO_2$ from \cite{Cun1}
with respect to representation theory.

A common thread of the constructions in \cite{Cun1}, \cite{CuLi1}, \cite{CuLi2} and \cite{Li} is
that the associated $C^*$-algebra $\fA[R]$ of the ring $R$ is a semigroup crossed
product of a commutative $C^*$-algebra $\cD$ by an action of the semidirect product semigroup $P_R=R\rtimes R^\times$
with $R^\times =R\setminus \{0\}$. Thus
by the dilation theory of \cite{Laca} the algebra $\fA[R]$ is Morita equivalent to an ordinary crossed product $\fA(R)$  by the group
$P_{Q(R)}=Q(R)\rtimes Q(R)^\times$, where $Q(R)$ is the quotient field of $R$. Here $P_{Q(R)}$ acts on a
commutative $C^*$-algebra with the locally compact spectrum $\mathscr{R}$ equal to a completion of $\Spec(\cD)$.
In the case of the ring of integers in a global field (a number field
or a function field), $\mathscr{R}$ is the ring of finite adeles, and one of the crucial
results in \cite{CuLi2} is a duality theorem which establishes a Morita equivalence between $\fA(R)$ and the crossed product of
 $P_{Q(R)}$ acting by affine transformations on the algebra of continuous functions on the space of infinite adeles.

Here we study the $C^*$-algebra, which we call $\cQ_2$, associated with the semidirect product
semigroup of the additive group $\Zz$ acted upon by multiplication with non-negative powers of $2$. We show that $\cQ_2$ admits a
crossed product structure mirroring the structure of $\cQ_\Nz$ at the prime $p=2$. For this reason, we refer to this $C^*$-algebra as
\emph{the $2$-adic ring $C^*$-algebra of the integers}. We establish a duality theorem in the form of
a Morita equivalence between the crossed products arising from the affine actions of the group $\Zz[\halb]\rtimes \spkl{2}$,
the semidirect product arising from the multiplicative action of the subgroup of $\Qz \reg$ generated by $2$ on $\Zz[\halb]$,
on $C_0(\Qz_2)$ and $C_0(\Rz)$. Both of these crossed products are Morita equivalent to $\cQ_2$. Appealing to Zhang's dichotomy we can argue that the
$C^*$-algebras $C_0(\Rz)\rtimes \Zz[\halb]\rtimes \spkl{2}$ and $C_0(\Qz_2)\rtimes \Zz[\halb]\rtimes \spkl{2}$ are isomorphic, although
we do not have an explicit isomorphism at hand, see Remark~\ref{Zhang-dichotomy}.

However, we obtain more than just the analogues of the structural results on $\cQ_\Nz$. The first additional information is
a complete characterisation of those representations of $\cO_2$ which extend to $\cQ_2$. The result motivates viewing
$\cQ_2$ as a symmetrized version of $\cO_2$ because it emphasizes a certain bond between the two generating isometries.  Further,
we work out a concrete description of the imprimitivity bimodule which implements the Morita equivalence of
$\cQ_2$ and the crossed product of $C_0(\Rz)$ by  $\Zz[\halb]\rtimes \spkl{2}$. We then show that via this bimodule, the
canonical representation of $\cQ_2$ on $\ell^2(\Zz)$ is induced to a representation of $C_0(\Rz)\rtimes \Zz[\halb]\rtimes \spkl{2}$
 which is unitarily equivalent to the canonical representation coming from multiplication, translation and dilation by $2$ on $\Ell^2(\Rz).$

The $C^*$-algebra $\cQ_2$ has appeared earlier in other contexts, see \cite{Hir}, or \cite{EanHR} and the list of references therein.
We refer to Section~\ref{section-is} for more details.
Moreover, we point out that our new approach to this algebra presented in this paper not only provides alternative proofs of established results;
it also reveals a close connection to ring $C^*$-algebras and leads to new insights on the representation theory of this algebra.

We also remark that we have chosen the particular number 2 to be concrete. We
can replace 2 by an arbitrary prime number or even an arbitrary positive
integer $n$ in the definition of $\cQ_2$, and then the constructions as well
as the results and their proofs will carry over
(with appropriate modifications).

The first named author thanks Joachim Cuntz, Siegfried Echterhoff and the group in noncommutative geometry at the Westf{\"a}lische Wilhelms-Universit{\"a}t
M{\"u}nster for their hospitality during a sabbatical stay in October 2009. The second named author thanks the operator algebra group in Oslo for a nice visit at the University of Oslo. We thank Joachim Cuntz for suggesting the problem and for useful discussions. We thank the anonymous referee for suggesting changes that resulted in a clearer description of the relation between $\cQ_2$ and Elliott's  classification  program for $C^*$-algebras.

\section{Preliminaries and notation}

Recall from \cite{Cun2} that the ring $C^*$-algebra of the integers $\cQ_{\Zz}$ is the universal $C^*$-algebra
generated by isometries $s_m$ with $m \in \Zz^\times$ and a unitary
$u$ satisfying the relations $s_k s_m = s_{km}$, $s_m u^n=u^{mn} s_m$ and $\sum_{l=0}^{m-1} u^l s_m s_m^* u^{-l} = 1$ for $k,m \in \Zz^\times$ and $n \in \Zz$.
By \cite[Theorem~7.1]{Cun2}, $\cQ_{\Zz}$ is simple and purely infinite.

Further, it was shown in \cite{CuLi1} and \cite{CuLi2} that $\cQ_{\Zz}$ is a semigroup crossed product of $C(\widehat{\Zz})$ by $P_\Zz$,
where $\widehat{\Zz}$ is the compact ring of finite integral adeles and $P_\Zz=\Zz\rtimes \Zz^\times$. By Laca's dilation theory
(see \cite{Laca}), $\cQ_{\Zz}$ is  Morita equivalent to the group crossed product $C_0(\Az_f)\rtimes P_{\Qz}$, where $\Az_f$
is the ring of finite adeles and $P_{\Qz}$ acts by affine transformations, see \cite[Remark~2]{CuLi1}.

Note that the crossed product description and the property of $\fA[R]$ being simple and purely infinite
(where $\fA[\Zz]$ is $\cQ_\Zz$) are  valid for a large class of integral domains $R$,
cf. \cite[Theorem~1, Theorem~2 and Remark~2]{CuLi1},
and in fact for a large class of rings $R$, see \cite[Proposition~4.1 and \S~5.2]{Li}.

The duality theorem in the case of $\Zz$ in its most general form is \cite[Corollary~3.10]{CuLi2}, and
taking the subgroup $\Gamma$ to be $\Qz^\times$ gives a Morita equivalence between $C_0(\Rz)\rtimes P_{\Qz}$
and $C_0(\Az_f)\rtimes P_{\Qz}$, where both actions of $P_{\Qz}$ are by (inverse) affine transformations.

We are interested in the case of the prime $2$ acting on $\Zz$.
Let $[2\rangle$ denote the multiplicative subsemigroup generated by the number $2$ in $\Zz \reg = \Zz \setminus \gekl{0}$,
i.e. $[2\rangle = \gekl{1,2^1,2^2,2^3,\dotsc}$,
and let $\spkl{2}$ denote the multiplicative subgroup $\menge{2^i}{i \in \Zz}$ of $\Zz [\halb] \reg$,
where $\Zz [\halb] \defeq [2\rangle^{-1} \Zz = \menge{\tfrac{l}{2^i}}{l \in \Zz,i \in \Nz_0}$.

The ring of 2-adic integers $\Zz_2$ is identified with $\plim_i \gekl{\Zz / 2^i \Zz}$
where the structure maps are the canonical projections $\Zz / 2^{i+1} \Zz \onto \Zz / 2^i \Zz$ given by reduction $\operatorname{mod} 2^i$.
As usual, $\Qz_2$ denotes the field of 2-adic numbers, i.e. the quotient field of $\Zz_2$.
We use symbols such as $\zf$ for elements in $\Zz_2$ and $\xf$, $\yf$ for elements in $\Qz_2$.

The semidirect product of the additive group $\Zz$ by the multiplicative semigroup $[2\rangle$ is the semigroup $\Zz\rtimes [2\rangle$ with operation $(l, 2^i)(n, 2^j)=(l+2^i n, 2^{i+j})$.  Then $\Zz\rtimes [2\rangle$ is a right-reversible Ore semigroup because
\bgloz
  (-2^j l,2^j)(l,2^i) = (0,2^{j+i}) = (0,2^{i+j}) = (-2^i n,2^i)(n,2^j).
\egloz
The corresponding enveloping group of left-quotients is $\Zz [\halb] \rtimes \spkl{2}$ in which the operation is given by $(b,a)(d,c) = (b+ad,ac)$ for $b,d \in \Zz [\halb]$ and $a,c \in \spkl{2}$. We shall use notation like $(b,a)$, $(d,c)$ for elements in $\Zz [\halb] \rtimes \spkl{2}$. We will use the fact that the
isomorphism $\spkl{2} \ni a \ma a^{-1} \in \spkl{2}$ implements an isomorphism $(b,a) \ma (b,a^{-1})$
between the semidirect product $\Zz[\halb] \rtimes \spkl{2}$ with respect to the ordinary multiplicative action
and the semidirect product with respect to the action by inverse multiplication.

Let $P$ be a discrete semigroup with an identity element, and let $P$ act by endomorphisms on a unital $C^*$-algebra $D$ via an action $\alpha: P \ri \End(D)$ that preserves the identity elements. The corresponding semigroup crossed product $D \rte_{\alpha} P$ is the
universal $C^*$-algebra for covariant pairs $(\pi, W)$ consisting  of a representation $\pi$ of $D$ on a Hilbert space $H$ and a homomorphism $W$ from $P$ into the semigroup of isometries on $H$ such that the \emph{covariance condition} $\pi(\alpha_p(a))=W_p \pi(a) W_p^*$ is satisfied for all $a\in D$ and $p\in P$,
see e.g. \cite{LacRae}. We write $(\iota_D, w)$ for the universal covariant representation, and we let $\menge{w_p}{p \in P}$ be the isometries in $D \rte_{\alpha} P$ which implement  $\alpha$.

We let $\e(t) \defeq \exp(2 \pi i t)$ for $t\in \Rz$.

\section{The 2-adic ring $C^*$-algebra of the integers}
\label{Q2}

\bdefin
Let $\cQ_2$ denote the universal unital $C^*$-algebra  generated by a unitary $u$ and an isometry $s_2$ subject to the relations
\bgloz
  (I) \ s_2 u = u^2 s_2,
\egloz
\bgloz
  (II) \ s_2 s_2^* + u s_2 s_2^* u^{-1} = 1.
\egloz
\edefin

If we compare the definition of $\cQ_2$ with the definition of the ring $C^*$-algebra $\fA[R]$ of a
ring $R$ from \cite{CuLi1} or \cite{Li} it is clear that $\cQ_2$ is not of the form $\fA[R]$ for any ring, at least not directly.
However, it will become clear in Section~\ref{croprodes} that $\cQ_2$ is associated with the semigroup $\Zz\rtimes [2\rangle$. Thus
 $\cQ_2$ can be viewed as a quotient of the $C^*$-algebra of a semigroup $R \rtimes H$ (in the sense of \cite[Definition~11]{Li}), where $H$ is a multiplicatively closed
 subset of $R^\times$.

There is a canonical representation of $\cQ_2$ on $\ell^2(\Zz)$: Let $\menge{\varepsilon_n}{n \in \Zz}$ be the canonical orthonormal basis of $\ell^2(\Zz)$. Consider the unitary $U$ defined by $U \varepsilon_n = \varepsilon_{n+1}$ and the isometry $S_2$ given by $S_2 \varepsilon_n = \varepsilon_{2n}$. It is easy to check that relation (I) is satisfied for $U$ and $S_2$. To check relation (II), we first observe that $S_2^* \varepsilon_n = \1z_{2 \Zz}(n) \varepsilon_{\tfrac{n}{2}}$ where $\1z_{2 \Zz}$ is the characteristic function of $2 \Zz$, defined on $\Zz$. It follows that $S_2 S_2^* \varepsilon_n = \1z_{2 \Zz}(n) \varepsilon_n$ and $U S_2 S_2^* U^* \varepsilon_n = \1z_{1+2 \Zz}(n) \varepsilon_n$. Here, $\1z_{1+2 \Zz}$ is the characteristic function of the coset $1+2 \Zz \subseteq \Zz$. Thus, relation (II) is satisfied as it corresponds to the decomposition $\Zz = (2 \Zz) \dotcup (1+2 \Zz)$. By the universal property of $\cQ_2$ there is a representation of $\cQ_2$ on $\ell^2(\Zz)$ which sends $u$ to $U$ and $s_2$ to $S_2$. We let $\lambda_2$ be this representation.

Moreover, there is a canonical homomorphism from $\cQ_2$ into the ring $C^*$-algebra of the integers. This follows from the universal property of $\cQ_2$ since the elements $u$ and $s_2$ in $\cQ_{\Zz}$  satisfy relations (I) and (II). The above representation $\lambda_2$ of $\cQ_2$ is the composition of the canonical homomorphism $\cQ_2 \ri \cQ_{\Zz}$ and the left regular representation of $\cQ_{\Zz}$ on $\ell^2(\Zz)$.

\section{The inner structure}
\label{section-is}

In the sequel we study the inner structure of $\cQ_2$.  Though $\cQ_2$ is not an example of a ring $C^*$-algebra in the sense of \cite{Li}, it is purely infinite and simple. Recall that this means that for all non-zero $x$ in $\cQ_2$ there exist $y$, $z$ in $\cQ_2$ such that $yxz=1$. To establish this claim, we can follow the strategy of proof from \cite[Theorem~1]{CuLi1} (see also the Section \an{Ring $C^*$-algebras for commutative rings} in \cite{Li}). This is indicated below. Alternatively, we can also express $\cQ_2$ as a Cuntz-Pimsner algebra and invoke \cite[Proposition~4.2]{EanHR}, as shown in Remark~\ref{other_contexts}.

We begin with some observations that will be useful later on. For every $i$ in $\Nz_0$, we write $s_{2^i}$ for $(s_2)^i$. The range projection of $s_{2^i}$ is denoted by $e_{2^i}$, so $e_{2^i} = s_{2^i} s_{2^i}^*$. It follows from relation (I) that
\bgloz
  s_{2^i} u^l = u^{2^i l} s_{2^i} \fa i \in \Nz_0 \text{ and } l\in \Zz.
\egloz
Moreover, relations (I) and (II) imply that we have
\bgl
\label{sumproj}
  \sum_{l \in \Zz / 2^i \Zz} u^l e_{2^i} u^{-l} = 1
\egl
for each $i$ in $\Nz_0$, or more generally,
\bgl
\label{sumprojgen}
  \sum_{l \in 2^i \Zz / 2^j \Zz} u^l e_{2^j} u^{-l} = e_{2^i}
\egl
for all integers $j \geq i \geq 0$.

 We let $\cD$ be the $C^*$-subalgebra of $\cQ_2$ generated by $\menge{u^l e_{2^i} u^{-l}}{l \in \Zz, i \in \Nz_0}$. Either by integrating dual actions analogously to \cite[Proposition~1]{CuLi1}, or by using the crossed product description of $\cQ_2$ developed in the next section as in \cite{Li}, we obtain a faithful conditional
expectation $\Theta: \cQ_2 \ri \cD$ which is characterized by
\bgloz
  \Theta(s_{2^i}^* u^{-l} f u^{l'} s_{2^{i'}}) = \delta_{i,i'} \delta_{l,l'} s_{2^i}^* u^{-l} f u^l s_{2^i} \fa i,i' \in \Nz_0; l,l' \in \Zz; f \in \cD.
\egloz
One can then describe this expectation by well-chosen projections as in \cite[Proposition~2]{CuLi1} or in \cite{Li}, Section \an{Ring $C^*$-algebras for commutative rings}
and follow the argument from \cite[Theorem~1]{CuLi1}, or in \cite{Li}, Section \an{Ring $C^*$-algebras for commutative rings} to establish the following result (we
omit the details).
\btheo
\label{pis}
$\cQ_2$ is purely infinite and simple.
\etheo

As a consequence of  simplicity of $\cQ_2$, we deduce that every non-zero representation of $\cQ_2$ is faithful. In particular, $\lambda_2$ is a faithful representation of $\cQ_2$ on $\ell^2(\Zz)$. Moreover, the canonical homomorphism $\cQ_2 \ri \cQ_{\Zz}$ must be injective, so that $\cQ_2$ can be identified with the corresponding $C^*$-subalgebra of $\cQ_{\Zz}$.

In view of the classification program for C*-algebras, it is important to note that $\cQ_2$ is
nuclear and satisfies the UCT, see Remark~\ref{other_contexts} or
Remark~\ref{nuclear-UCT}. This means that $\cQ_2$ is a Kirchberg algebra which satisfies the UCT, hence it belongs to a class of C*-algebras which is classified by K-theory (see \cite{Ror}, Chapter~8). Moreover, the K-theory of $\cQ_2$ can be determined and the class of the unit $1_{\cQ_2}$ in K-theory vanishes, see \eqref{K(Q_2)}.

\bremark
\label{other_contexts}
The $C^*$-algebra $\cQ_2$ has appeared in other contexts, see \cite{Hir}, \cite{Ka}, or \cite{EanHR} and the
references given therein. In \cite{Hir}, the $C^*$-algebra
$\cE_{\alpha}$  associated with the
 endomorphism $\alpha = 2 \id_{\Zz}$ of $\Zz$ is generated by an isometry
and a unitary subject to relations which are precisely (I) and (II), and is simple by
\cite[Theorem~3.9]{Hir}.

In the context of Katsura's extensive work on
$C^*$-algebras associated to topological graphs,
the algebra $\cQ_2$ is the same as $\mathcal{O}(E_{2,1})$ from \cite[Example A.6]{Ka} and hence is a Kirchberg algebra. Moreover,
\begin{equation}
\label{K(Q_2)}
K_0(\cQ_2)\cong \Zz \text{ and } K_1(\cQ_2)\cong\Zz,
\end{equation}
and the class of the unit $1_{\cQ_2}$ vanishes in $K_0(\cQ_2)$. This can be derived from relation (II) as follows:
\bgloz
  \eckl{1} \overset{(II)}{=} \eckl{s_2 s_2^* + u s_2 s_2^* u^{-1}} = \eckl{s_2 s_2^*} + \eckl{u s_2 s_2^* u^{-1}} = \eckl{1} + \eckl{1}.
\egloz
Here $\eckl{\cdot}$ stands for $K_0$-class. Alternatively, one can also use the unital embedding $\cO_2 \into \cQ_2$ (see the next section) and the induced map on $K_0$ to deduce that the class of $1_{\cQ_2}$ vanishes in $K_0(\cQ_2)$.

It was explained in the introduction to \cite{EanHR} that the crossed product of an Exel system
$(C(\Tz), \sigma, L, \Nz)$  gives rise to an algebra of the sort studied in \cite[Appendix A]{Ka}. Here we prove directly that the Exel crossed product $C(\Tz)\rtimes_{\sigma, L}\Nz$
is isomorphic to $\cQ_2$, hence to the afore-mentioned $\mathcal{O}(E_{2,1})$. Towards this, we recall that $\sigma$ is the endomorphism of $C(\Tz)$  given through the
covering map $z\to z^2$ on $\Tz$ by $\sigma(f)(z)=f(z^2)$. Then $M=C(\Tz)$ is a right
$C(\Tz)$-module with the action
$x\cdot f=x\sigma(f)$ for $x, f\in C(\Tz)$. The map
$L:C(\Tz)\to C(\Tz)$, $L(f)(z)=\frac 12\sum_{w^2=z}f(w)$ satisfies $L(f\sigma(g))=L(f)g$ for $f,g\in C(\Tz)$, hence is a transfer operator for $\sigma$ in the sense of Exel, see \cite{Exel}.
Then $\langle x,y\rangle:= L(x^*y)$ gives rise to a $C(\Tz)$-valued inner product on $M$.
It was shown in \cite[Lemma~3.3]{LarRae2} that
$C(\Tz)$ is complete in the norm $\Vert x\Vert=\langle x,x\rangle^{1/2}$, hence $M$ is a right Hilbert $C(\Tz)$-module. There is a left action
implemented by the homomorphism $\phi:C(\Tz)\to \cL(M)$  given by
$(\phi(f)x)(z)=f(z) x(z)$. By \cite[\S~2.2.]{EanHR}, the Exel crossed product $C(\Tz)
\rtimes_{\sigma, L}\Nz$  associated
with $(C(\Tz), \sigma, L, \Nz)$ is the Cuntz-Pimsner algebra $\mathcal{O}_M$.

To see that $C(\Tz)\rtimes_{\sigma, L}\Nz$ is isomorphic to $\cQ_2$ let $\id$ be the identity
function $z\mapsto z$  in $C(\Tz)$ and in $M$, and write $\1z$ for the
constant function $1$, seen as an element of $M$. The universal Toeplitz representation of $M$ is a pair consisting of a linear map $\iota_M:M\to \mathcal{T}_M$ and a
$^*$-homomorphism $\iota_A:C(\Tz)\to \mathcal{T}_M$ which are compatible with the module operations and the inner-product, cf. \cite{Pims}.

Now set $\tilde{s}=\iota_M(\1z)$ and $\tilde{u}=\iota_A(\id)$. Since $\id$ is a unitary in $C(\Tz)$, so is $\tilde{u}$ in
 $\mathcal{T}_M$, and since $L(1)=1$, the relation $\iota_A(\langle \1z,\1z \rangle)=\iota_M(\1z)^*\iota_M(\1z)$ shows that $\tilde{s}$
 is an isometry. Notice that the right action $\1z\cdot \id$ gives the same element in $M$ as the left action $\id^2\cdot \1z$, so
 $\tilde{s}\tilde{u}=\tilde{u}^2\tilde{s}$.

Denote by $\iota^{(1)}$ the homomorphism from $\cK(M)$ to $\cT_M$ such that $\iota^{(1)}(\theta_{x,y})$ equals $\iota_M(x)\iota_M(y)^*$
 when $x,y \in M$ and $\iota_A(x)\iota_A(y)^*$ when $x,y\in A$. In $\cL(M)$ we have $\phi(1)=\theta_{\1z, \1z}+\theta_{\id \cdot \1z, \id \cdot \1z}$. In particular
 $\mathcal{K}(M)=\cL(M)$, and Cuntz-Pimsner covariance amounts  to $\iota^{(1)}(\phi(f))=\iota_A(f)$ for all $f\in C(\Tz)$. It suffices
 to see what happens at $f=1$, and this amounts to the relation
 $1=\tilde{s}\tilde{s}^*+\tilde{u}\tilde{s}\tilde{s}^*\tilde{u}^*$. Thus the universal property of $\cQ_2$
 gives a homomorphism onto $\cO_M$ which is injective by simplicity of $\cQ_2$.
 \eremark

\section{Extending representations from $\cO_2$ to $\cQ_2$}

In this section we initiate the study of representations of $\cQ_2$. Recall first that $\cO_2$ is the universal $C^*$-algebra generated by two isometries $s_0$ and $s_1$ satisfying the relation $s_0s_0^*+s_1s_1^*=1$, see \cite{Cun1}. The next proposition says that $\cQ_2$ can be thought of as a symmetrized version of $\cO_2$, in the sense that the two generating isometries have to be unitarily equivalent, at least in concrete representations. The following is the first main result of this paper.

\bprop
\label{fromO2toQ2}
Let $S_0$ and $S_1$ be isometries on a Hilbert space $\cH$ which give rise to a representation of $\cO_2$, i.e. we have
\bgl
\label{S0S1}
  S_0 S_0^* + S_1 S_1^* = I.
\egl
Then there exists a representation $\pi$ of $\cQ_2$ on $\cH$ with
\bgloz
  \pi(s_2) = S_0 \text{ and } \pi(u s_2) = S_1
\egloz
if and only if the unitary parts in the Wold decompositions of $S_0$ and $S_1$ are unitarily equivalent.
\eprop

Recall that the unitary part of an isometry $S$ on some Hilbert space is given by the restriction of $S$ to the $S$-invariant subspace $\bigcap_{n=1}^{\infty} \img(S^n)$.

\bproof
By definition of $\cQ_2$, a representation $\pi$ of $\cQ_2$ with
\bgloz
  \pi(s_2) = S_0 \text{ and } \pi(u s_2) = S_1
\egloz
exists if and only if there is a unitary $U$ on $\cH$ with the properties
\bgln
\label{U1}
  U S_0 &=& S_1, \text{ and} \\
\label{U2}
  S_0 U &=& U^2 S_0.
\egln

If there is a unitary $U$ satisfying \eqref{U1} and \eqref{U2}, then $S_0 U = U U S_0 = U S_1$. Thus $S_0$ and $S_1$ themselves are unitarily equivalent. In particular, their unitary parts must be unitarily equivalent, and the \an{only if}-part follows.

It remains to prove the \an{if}-part. By assumption, there is a unitary operator
\bgloz
  \bigcap_{n=1}^{\infty} \img({S_1}^n) \overset{\cong}{\lori} \bigcap_{n=1}^{\infty} \img({S_0}^n)
\egloz
transforming $S_1$ into $S_0$. Extending this unitary operator by $0$ on the orthogonal complements, we obtain a partial isometry $W$ on $\cH$ with initial space $\bigcap_{n=1}^{\infty} \img({S_1}^n)$ and final space $\bigcap_{n=1}^{\infty} \img({S_0}^n)$. Moreover,
\bgl
\label{WS=SW}
  W S_1 = S_0 W.
\egl
Let us produce a partial isometry $V$ with initial space $\rukl{\bigcap_{n=1}^{\infty} \img({S_1}^n)}^{\perp}$ and final space $\rukl{\bigcap_{n=1}^{\infty} \img({S_0}^n)}^{\perp}$ so that $U \defeq V + W$ is a unitary with the desired properties.

Consider for every non-negative integer $n$ the bounded operator
\bgl
\label{Vn}
  V_n \defeq \sum_{i=0}^{n} S_0^i S_1 S_0^* {S_1^*}^i.
\egl
We claim that the sequence $V_n$ is strongly convergent. Then the strong limit will be the operator $V$ we are looking for.

For every non-negative integer $n$, we have
\bgln
\label{S1S1*}
  S_1 S_1^* &=& \sum_{j=1}^{n} \rukl{S_1^j {S_1^*}^j - S_1^{j+1} {S_1^*}^{j+1}} + S_1^{n+1} {S_1^*}^{n+1} \\
  &=& \sum_{j=1}^{n} S_1^j (I - S_1 S_1^*) {S_1^*}^j + S_1^{n+1} {S_1^*}^{n+1} \nonumber \\
  &=& \sum_{j=1}^{n} S_1^j S_0 S_0^* {S_1^*}^j + S_1^{n+1} {S_1^*}^{n+1} \text{ by } \eqref{S0S1}. \nonumber
\egln
Thus for every $n$ in $\Nz_0$
\bgl
\label{1=sum}
  I = S_0 S_0^* + S_1 S_1^* = \sum_{j=0}^{n} S_1^j S_0 S_0^* {S_1^*}^j + S_1^{n+1} {S_1^*}^{n+1}.
\egl
It then follows that
\bgln
\label{Vn*Vn}
  V_n^* V_n &=& \sum_{i,j=0}^{n} S_1^j S_0 S_1^* {S_0^*}^j S_0^i S_1 S_0^* {S_1^*}^i \\
  &=& \sum_{j=0}^n S_1^j S_0 S_0^* {S_1^*}^j = I - S_1^{n+1} {S_1^*}^{n+1} \text{ by } \eqref{1=sum}. \nonumber
\egln
In particular, $V_n^* V_n \leq 1$, so that
\bgl
\label{normVn}
  \norm{V_n} \leq 1 \fa n \in \Nz_0.
\egl
By \eqref{S0S1}, we know that for every vector $\xi$ in $\cH$, the vectors
\bgloz
  S_0^i S_1 S_0^* {S_1^*}^i \xi, \ i \in \Nz_0
\egloz
are pairwise orthogonal. Thus we have for $n \geq m$:
\bgln
\label{Vn-Vm}
  && \norm{V_n \xi}^2 = \sum_{i=0}^n \norm{S_0^i S_1 S_0^* {S_1^*}^i \xi}^2 \\
  &=& \norm{V_m \xi}^2 + \sum_{i=m+1}^n \norm{S_0^i S_1 S_0^* {S_1^*}^i \xi}^2 = \norm{V_m \xi}^2 + \norm{V_n \xi - V_m \xi}^2. \nonumber
\egln
So $(\norm{V_n \xi}^2)_n$ is a bounded (see \eqref{normVn}) and monotonically increasing sequence. Thus it converges, which means that $(\norm{V_n \xi}^2)_n$ is a Cauchy sequence. By \eqref{Vn-Vm}, this implies that $(V_n \xi)_n$ is a Cauchy sequence in norm. Hence for every vector $\xi$ in $\cH$, the sequence $(V_n \xi)_n$ converges. If we define $V: \cH \ri \cH$ by $V \xi \defeq \lim_{n \ri \infty} V_n \xi$ for all $\xi$ in $\cH$, then we obtain a linear operator which by \eqref{normVn} satisfies
\bgloz
  \norm{V \xi} = \lim_{n \ri \infty} \norm{V_n \xi} {\leq} \norm{\xi}
\egloz
for all $\xi$ in $\cH$. Thus $V$ is bounded. In other words, the sequence $(V_n)_n$ does indeed converge strongly and $V$ is the strong limit of $(V_n)_n$.

By construction, $V_n^*$ is of the same form as $V_n$, but with the roles of $S_0$ and $S_1$ interchanged. Thus by the same argument, but with reversed roles for $S_0$ and $S_1$, we obtain that the sequence $(V_n^*)_n$ converges strongly to some operator. The strong limit of $(V_n^*)_n$ must then coincide with $V^*$.

Using sequential continuity of multiplication with respect to the strong operator topology, we deduce from \eqref{Vn*Vn} that
\bgloz
  V^* V = \lim_{n \ri \infty} V_n^* V_n = \lim_{n \ri \infty} (I - S_1^{n+1} {S_1^*}^{n+1})
\egloz
is the orthogonal projection onto $\rukl{\bigcap_{n=1}^{\infty} \img({S_1}^n)}^{\perp}$. Analogously, we obtain that $V V^*$ is the orthogonal projection onto $\rukl{\bigcap_{n=1}^{\infty} \img({S_0}^n)}^{\perp}$. Thus $V$ is indeed a partial isometry with initial space $\rukl{\bigcap_{n=1}^{\infty} \img({S_1}^n)}^{\perp}$ and final space $\rukl{\bigcap_{n=1}^{\infty} \img({S_0}^n)}^{\perp}$. Set $U \defeq V + W$. By construction, it is clear that $U$ is a unitary.

Now, $W$ is $0$ on $\rukl{\bigcap_{n=1}^{\infty} \img({S_1}^n)}^{\perp} \supseteq \img(S_1)^{\perp} = \img(S_0)$ by construction. Thus, using sequential continuity of multiplication in the strong operator topology, we infer that
\bgln
\label{U1v}
  U S_0 &=& V S_0 + W S_0 = V S_0 = \lim_{n \ri \infty} V_n S_0 \\
  &=& \lim_{n \ri \infty} \sum_{i=0}^n S_0^i S_1 S_0^* {S_1^*}^i S_0 = \lim_{n \ri \infty} S_1 S_0^* S_0 = S_1 \nonumber
\egln
and
\bglnoz
  S_0 U &=& S_0 V + S_0 W = \rukl{\lim_{n \ri \infty} S_0 V_n} + W S_1 \text{ by } \eqref{WS=SW}\\
  &=& \lim_{n \ri \infty} \rukl{\sum_{i=0}^n S_0^{i+1} S_1 S_0^* {S_1^*}^i} + W S_1 \\
  &=& \lim_{n \ri \infty}\rukl{\sum_{i=0}^{n+1} S_0^i S_1 S_0^* {S_1^*}^i} S_1 + W S_1 \\
  &=& (\lim_{n \ri \infty} V_{n+1} + W) S_1 = U S_1 = U^2 S_0 \text{ by } \eqref{U1v}.
\eglnoz
Therefore the operators $U$ and $S_0$ give rise to a representation of $\cQ_2$ with the properties \eqref{U1} and \eqref{U2}. This completes the proof of the proposition.
\eproof

\bremark
\label{rmk-unique}
It can be additionally deduced from our arguments that the restriction of the unitary to $\rukl{\bigcap_{n=1}^{\infty} \img({S_1}^n)}^{\perp}$ is uniquely determined: On this subspace, it will always coincide with the operator $V$ we constructed in the proof.
\eremark

We may rephrase our first main result in the following way: $\cO_2$ sits as a $C^*$-subalgebra in $\cQ_2$ via the homomorphism sending $s_0$ to $s_2$ and $s_1$ to $u s_2$. Proposition~\ref{fromO2toQ2} tells us when precisely a representation of this $C^*$-subalgebra extends to a representation of the whole $C^*$-algebra.
Remark~\ref{rmk-unique} explains to which extent this extension is unique.

Our extension result can be applied in the following context: Bratteli and Jorgensen studied connections between wavelets and certain representations of Cuntz algebras in
\cite{BraJor}. Proposition~\ref{fromO2toQ2} tells us precisely which of the representations of $\cO_2$ can be extended to $\cQ_2$. For example, we obtain that the representations of $\cO_2$ given by the isometries $S_0$ and $S_1$ in \cite[Equation~(1)]{LarRae1} extend to $\cQ_2$. Namely, \cite[Theorem~3.1]{BraJor} implies that both isometries $S_0$ and $S_1$ are pure, so they both have zero unitary part.

\section{$\cQ_2$ as a semigroup crossed product}

\label{croprodes}

In the following we represent $\cQ_2$ as a semigroup crossed product. This observation builds the bridge to the 2-adic integers and the 2-adic numbers. This then justifies why we can think of $\cQ_2$ as the 2-adic version of the ring $C^*$-algebra of the integers.

Recall that  $ \cD$ is the closed span of projections $u^l e_{2^i} u^{-l}$ for ${l \in \Zz, i \in \Nz_0}$. If we let $\Ad(u^l s_{2^i})(f) = u^l s_{2^i} f (u^l s_{2^i})^*$ for $f\in \cD$, it follows from relation (I) that $\Ad(u^l s_{2^i}) (\cD) \subseteq \cD$ for all $i$ in $\Nz_0$ and $l$ in $\Zz$. This defines an action $\alpha$ of the semigroup $\Zz \rtimes [ 2 \rangle$ on $\cD$ via
\bgloz
  \Zz \rtimes [ 2 \rangle \ni (l,2^i) \ma \Ad(u^l s_2^i) \in \End(\cD).
\egloz

\bprop
\label{Q2semigpcp}
The map $\cQ_2 \rightarrow\cD \rte_{\alpha} (\Zz \rtimes [ 2 \rangle)$ given by
\bgl
\label{Q2semi''}
  u^l s_2^i \ma w_{(l,2^i)}
\egl
for $l\in \Zz$ and $i\in \Nz_0$ is an isomorphism.
\eprop

\bproof
We compare the universal properties. The isometry $w_{(0,2)}$ and the unitary $w_{(1,1)}$ satisfy relation (I) by covariance, and relation (II) follows from
\bgloz
  w_{(0,2)} w_{(0,2)}^* + w_{(1,2)} w_{(1,2)}^* = \alpha_{(0,2)}(1) + \alpha_{(1,2)}(1) = e_2 + u e_2 u^{-1} = 1.
\egloz
By  the universal property of $\cQ_2$, the map \eqref{Q2semi''} defines a homomorphism $\cQ_2 \ri \cD \rte_{\alpha} (\Zz \rtimes [ 2 \rangle)$, and this map is injective by simplicity of $\cQ_2$.

We claim that $\cD \rte_{\alpha} (\Zz \rtimes [ 2 \rangle)$ is generated as a $C^*$-algebra by $w_{(0,2)}$ and $w_{(1,1)}$. As a semigroup crossed product, $\cD \rte_{\alpha} (\Zz \rtimes [ 2 \rangle)$ is the closed span of monomials $w_{(l, 2^i)}^* \iota_{\cD}(f) w_{(n, 2^j)}$ with $f \in \cD$ and $(l, 2^i), (n, 2^j) \in \Zz \rtimes [ 2 \rangle$ (this holds true by \cite{Laca}, Remark~1.3.1 and because $\Zz \rtimes [2 \rangle$ is an Ore semigroup). By the covariance relation,
\bgloz
\iota_{\cD}(u^l e_{2^i} u^{-l})=\iota_{\cD}(\alpha_{(l, 2^i)}(1))=w_{(l, 2^i)} w_{(l, 2^i)}^*,
\egloz
and since elements of the form $u^l e_{2^i} u^{-l}$ span $\cD$ and $w_{(l, 2^i)} = w_{(1,1)}^l w_{(0,2)}^i$, the claim follows. Hence the map in \eqref{Q2semi''} is an isomorphism.
\eproof

Next we show that $\cD$ is commutative and that the Gelfand transform on $\cD$ gives rise to an isomorphism of $\cD \rte_{\alpha} (\Zz \rtimes [ 2 \rangle)$ with $C(\Zz_2) \rte_{\alpha^{\operatorname{aff}}} (\Zz \rtimes [ 2 \rangle)$ where
\bgloz
  \alpha^{\operatorname{aff}}_{(l,2^i)} (f)(\zf) =
  \begin{cases}
    f(2^{-i}(\zf -l)), &\text{ if } \zf \in l+ 2^i \Zz_2 \\
    0, &\text{ otherwise}.
  \end{cases}
\egloz

\bprop
\label{DCZ2}
The $C^*$-algebra $\cD$ is commutative, and the Gelfand transform has range isomorphic to $C(\Zz_2)$ and maps $u^l e_{2^i} u^{-l}$ to the characteristic function $\1z_{l+2^i \Zz_2}$ of the subset $l+2^i \Zz_2$ of $\Zz_2$. Further, the Gelfand transform is equivariant for the actions $\alpha$ and $\alpha^{\operatorname{aff}}$.
\eprop

\bproof
We showed in \eqref{sumprojgen} that the projection $e_{2^i}$ can be written as a finite sum of the projections $\menge{u^l e_{2^j} u^{-l}}{l \in \Zz}$ for all integers $j \geq i \geq 0$. Let $\cD_i = C^*(\menge{u^l e_{2^i} u^{-l}}{l \in \Zz})$. There is therefore an inclusion $\iota_{i,i+1}: \cD_i \into \cD_{i+1}$, and hence an identification
\bgl
\label{D=ilim}
  \cD \cong \ilim_{i \in \Nz_0} \gekl{\cD_i; \iota_{i,i+1}}.
\egl
For fixed $i$, \eqref{sumproj} shows that the projections of the form $u^l e_{2^i} u^{-l}$ are pairwise orthogonal. Thus $\cD_i$ is a direct sum of $2^i$ copies of $\Cz$, hence commutative for every $i \in \Nz_0$. Therefore \eqref{D=ilim} implies that $\cD$ is commutative and the spectrum of $\cD$ is the inverse limit
\bgloz
  \Spec(\cD) \cong \plim_i \gekl{\Spec(\cD_i)); (\iota_{i,i+1})^*}
\egloz
where $\iota_{i,i+1}^*(\chi) = \chi \circ \iota_{i,i+1}$. Let $\iota_i: \cD_i \into \cD$ be the canonical injections. Then $(\iota_i)^*: \Spec(\cD) \ri \Spec(\cD_i)$ is the canonical restriction map given by $(\iota_i)^*(\chi) = \chi \circ \iota_i$. We can identify $\Zz / 2^i \Zz$ with $\Spec(\cD_i)$ via the map
\bgloz
  \sigma_i: l+2^i \Zz \ma \eckl{u^n e_{2^i} u^{-n} \ma \delta_{l+2^i \Zz, n+2^i \Zz}}.
\egloz
The canonical projections $p_{i+1,i}: \Zz / 2^{i+1} \Zz \onto \Zz / 2^i \Zz$ give rise to an identification of $\Zz_2$ as the profinite limit $\plim_i \gekl{\Zz / 2^i \Zz; p_{i+1,i}}$. Let $p_i: \Zz_2 \to \Zz / 2^i \Zz$ denote the canonical projection. Since $\sigma_i \circ p_{i+1,i}=(\iota_{i,i+1})^* \circ \sigma_{i+1}$, the maps
$\sigma_i$ induce an isomorphism $\sigma: \Zz_2 \to \Spec(\cD)$ such that $(\iota_i)^* \circ \sigma = \sigma_i \circ p_i$ for all $i \in \Nz_0$. The Gelfand map $\Gamma_{\cD}$ is an isomorphism $\cD \to C(\Spec(\cD))$. We claim that the composition of $\sigma^*: C(\Spec(\cD)) \to C(\Zz_2)$ given by $\sigma^*(f)= f \circ \sigma$ and $\Gamma_{\cD}$ is the isomorphism $\cD \cong C(\Zz_2)$ with
\bgl
\label{DC(Z2)}
  (\sigma^* \circ \Gamma_{\cD})(u^l e_{2^i} u^{-l})=\1z_{l+2^i \Zz_2}.
\egl
Indeed, for $\zf \in \Zz_2$ we have
\bglnoz
  && (\sigma^* \circ \Gamma_{\cD})(u^l e_{2^i} u^{-l})(\zf)=\Gamma_{\cD}(u^l e_{2^i} u^{-l})(\sigma(\zf))=\sigma(\zf)(u^l e_{2^i} u^{-l}) \\
  && = \sigma(\zf)(\iota_i({u^l e_{2^i} u^{-l}})) =((\iota_i)^* \circ \sigma)(\zf))(u^l e_{2^i} u^{-l}) \\
  && = \sigma_i(p_i(\zf))(u^l e_{2^i} u^{-l}) = \delta_{p_i(\zf),l+2^i \Zz} = \1z_{l+2^i\Zz_2}(\zf)
\eglnoz
where in the last equality we invoke the isomorphism $\Zz / 2^i \Zz \cong \Zz_2 / 2^i \Zz_2$ implemented by the maps
$l+2^i\Zz \ma l+2^i\Zz_2$, $p_i(\zf) \mafr \zf+2^i\Zz_2$. This implies \eqref{DC(Z2)}.

To see that $\sigma^* \circ \Gamma_{\cD}$ carries $\alpha$ into $\alpha^{\operatorname{aff}}$ we first compute that
\bgloz
  \alpha_{(n,2^j)}(u^l e_{2^i} u^{-l}) = u^{n+2^j l} e_{2^{j+i}} u^{-(n+2^j l)},
\egloz
and then apply $\sigma^* \circ \Gamma_{\cD}$ to obtain
\bgloz
  \1z_{n+2^j l + 2^{j+i} \Zz_2} = \1z_{n+2^j \Zz_2}(2^{-i}(\sqcup - l)) \1z_{l+2^i \Zz_2} = \alpha^{\operatorname{aff}}_{(n,2^j)}(\1z_{l+2^i \Zz_2}).
\egloz
\eproof

Proposition~\ref{Q2semigpcp} and Proposition~\ref{DCZ2} imply the following
\bcor
\label{Q2Z2}
There is an isomorphism $\cQ_2 \to C(\Zz_2) \rte_{\alpha^{\operatorname{aff}}} (\Zz \rtimes [ 2 \rangle)$ sending $s_2$ to $w_{(0,2)}$ and $u$ to $w_{(1,1)}$.
\ecor

Let $\beta^{\operatorname{aff}}$ be the action of $\Zz[\halb] \rtimes \spkl{2}$ on $C_0(\Qz_2)$ by affine transformations
\bgl
\label{aff-transf-Q2}
  (\beta^{\operatorname{aff}}_{(b,a)}f)(\xf)=f(a^{-1}(\xf -b))
\egl
for $(b,a) \in \Zz[\halb] \rtimes \spkl{2}$ and $f \in C_0(\Qz_2)$.

\blemma
\label{M'}
The triple $(C_0(\Qz_2),\Zz[\halb] \rtimes \spkl{2}, \beta^{\operatorname{aff}})$ is the minimal automorphic dilation of $(C(\Zz_2),\Zz \rtimes [ 2 \rangle, \alpha^{\operatorname{aff}})$. In particular, $C(\Zz_2) \rte_{\alpha^{\operatorname{aff}}} (\Zz \rtimes [ 2 \rangle)$ and $C_0(\Qz_2) \rtimes_{\beta^{\operatorname{aff}}} (\Zz[\halb] \rtimes \spkl{2})$ are Morita equivalent.
\elemma

\bproof
We show that $(C_0(\Qz_2), \Zz[\halb] \rtimes \spkl{2}, \beta^{\operatorname{aff}})$ and  the embedding $\iota$ of $C(\Zz_2)$ into $C_0(\Qz_2)$ given by $f \ma f \1z_{\Zz_2}$ satisfy the two conditions for the minimal automorphic dilation of $(C(\Zz_2), \Zz \rtimes [2 \rangle, \alpha^{\operatorname{aff}})$ from \cite[Theorem~2.1.1]{Laca}.

The first condition is satisfied because
\bgloz
  \iota(\alpha^{\operatorname{aff}}_{(l,2^i)}(f)) = f(2^{-i}(\sqcup - l)) \1z_{\Zz_2}(2^{-i}(\sqcup - l)) \1z_{\Zz_2}
  = (f \1z_{\Zz_2})(2^{-i}(\sqcup - l)) = \beta^{\operatorname{aff}}_{(l,2^i)}(\iota(f)).
\egloz
The second condition means that
\bgloz
  \bigcup_{(l,2^i) \in \Zz \rtimes [ 2 \rangle} (\beta^{\operatorname{aff}}_{(l,2^i)})^{-1}(\iota(C(\Zz_2))) \text{ is dense in } C_0(\Qz_2)
\egloz
and can be proven using the Stone-Weierstrass Theorem.

Morita equivalence now follows from \cite[Theorem~2.2.1]{Laca}.
\eproof

Combining Corollary \ref{Q2Z2} with Lemma \ref{M'}, we deduce
\bcor
\label{Q2MQ2}
$\cQ_2$ is Morita equivalent to $C_0(\Qz_2) \rtimes_{\beta^{\operatorname{aff}}} (\Zz[\halb] \rtimes \spkl{2})$.
\ecor

\section{The stabilization of $\cQ_2$}

In the previous section, we have seen a crossed product description of $\cQ_2$, first as a semigroup crossed product and then, after a dilation process, as an ordinary crossed product by a group.

In this section we study another dilation of $\cQ_2$. This dilation may seem to be the more canonical construction; however, it turns out that in the end,
we get the same $C^*$-algebra as in Corollary~\ref{Q2MQ2}.

Define the $C^*$-algebra $\bar \cQ_2$ as the inductive limit of the direct system
\bgloz
  \cQ_2 \overset{\Ad(s_2)}{\lori} \cQ_2 \overset{\Ad(s_2)}{\lori} \dotso
\egloz
where $\Ad(s_2)$ is the endomorphism $x\mapsto s_2 x s_2^*$ of $\cQ_2$.

We prove first that $\bar \cQ_2$ is isomorphic to $\cK(\ell^2(\Nz_0)) \otimes \cQ_2$, and
then we identify $\bar \cQ_2$ as a crossed product by an affine-type action.

\blemma
The map $\vartheta_{1,2}: \cQ_2 \lori M_2(\cQ_2)$ defined by $u \ma
  \rukl{
  \begin{smallmatrix}
  0 & u \\
  1 & 0
  \end{smallmatrix}
  }
$
and
$s_2 \ma
  \rukl{
  \begin{smallmatrix}
  s_2 & u s_2 \\
  0 & 0
  \end{smallmatrix}
  }
$ is an isomorphism.
\elemma
\bproof
The idea is that the mutually orthogonal projections $e_2$ and $u e_2 u^{-1}$ correspond in $2 \times 2$-matrices to the projections
$
  \rukl{
  \begin{smallmatrix}
  1 & 0 \\
  0 & 0
  \end{smallmatrix}
  }
$
and
$
  \rukl{
  \begin{smallmatrix}
  0 & 0 \\
  0 & 1
  \end{smallmatrix}
  }
$.
Moreover, these two projections $e_2$ and $u e_2 u^{-1}$ are Murray-von Neumann equivalent via $u e_2$, and this partial isometry corresponds in $2 \times 2$-matrices to
$
  \rukl{
  \begin{smallmatrix}
  0 & 0 \\
  1 & 0
  \end{smallmatrix}
  }
$,
the canonical partial isometry with support projection
$
  \rukl{
  \begin{smallmatrix}
  1 & 0 \\
  0 & 0
  \end{smallmatrix}
  }
$
and range projection
$
  \rukl{
  \begin{smallmatrix}
  0 & 0 \\
  0 & 1
  \end{smallmatrix}
  }
$.

The precise proof of the lemma goes as follows: An easy verification shows that relations (I) and (II) are valid for
$
  \rukl{
  \begin{smallmatrix}
  0 & u \\
  1 & 0
  \end{smallmatrix}
  }
$
and
$
  \rukl{
  \begin{smallmatrix}
  s_2 & u s_2 \\
  0 & 0
  \end{smallmatrix}
  }
$
in place of $u$ and $s_2$, respectively. Thus $\vartheta_{1,2}$ exists by the universal property of $\cQ_2$ and is injective by simplicity of $\cQ_2$. Let $\gekl{e_{i,j}}_{i,j=0,1}$ denote the matrix units of $M_2(\Cz)$. Then $\vartheta_{1,2}(e_2 u^{-1})=e_{0,1}$, $\vartheta_{1,2}(e_2)=e_{0,0}$, $\vartheta_{2,2}(ue_2u^{-1})=e_{1,1}$ and $\vartheta_{1,2}(ue_2)=e_{1,0}$. Thus $\vartheta_{1,2}$ is surjective.
\eproof

Next let $\vartheta_{2^i,2^{i+1}}: M_{2^i}(\cQ_2) \ri M_{2^{i+1}}(\cQ_2)$ be given by
\bgloz
  \vartheta_{2^i,2^{i+1}} \defeq \id_{2^i}\otimes \vartheta_{1,2}: M_{2^i}(\Cz) \otimes \cQ_2 {\lori}
  M_{2^i}(\Cz) \otimes M_2(\cQ_2).
\egloz
Thus $\vartheta_{2^i,2^{i+1}}$ acts as $\vartheta_{1,2}$ on every entry of $M_{2^i}(\cQ_2)$. For $i \geq 0$ let
\bgloz
  \vartheta_{2^i}
  \defeq \vartheta_{2^{i-1},2^i} \circ \vartheta_{2^{i-2},2^{i-1}} \circ \dotsb \circ \vartheta_{1,2} \text{ as a map }
  \vartheta_{2^i}:\cQ_2 {\overset{\cong}{\to}} M_{2^i}(\cQ_2).
\egloz
\blemma
Let $\rho_{2^i,2^{i+1}}: M_{2^i}(\cQ_2) \ri M_{2^{i+1}}(\cQ_2)$ be given by
$A \ma
  \rukl{
  \begin{smallmatrix}
  A & 0_{2^i} \\
  0_{2^i} & 0_{2^i}
  \end{smallmatrix}
  }
$. Then
\bgl
\label{thetaAd}
  \vartheta_{2^{i+1}} \circ \Ad(s_2) = \rho_{2^i,2^{i+1}} \circ \vartheta_{2^i}
\egl
for all $i$ in $\Nz_0$.
\elemma
\bproof
We proceed inductively on $i$. For $i=0$ we have $\vartheta_2 = \vartheta_{1,2}$ and $\vartheta_1 = \id$. Clearly $\vartheta_{1,2} \circ \Ad(s_2) = \rho_{1,2}$, and if \eqref{thetaAd} holds for $i-1$, then
\bglnoz
  &&\vartheta_{2^{i+1}} \circ \Ad(s_2) (x) = \vartheta_{2^i,2^{i+1}} \circ \vartheta_{2^i} \circ \Ad(s_2) (x) \\
  &=& \vartheta_{2^i,2^{i+1}} \circ \rho_{2^{i-1},2^i} \circ \vartheta_{2^{i-1}} (x) = \vartheta_{2^i,2^{i+1}}
  \rukl{
  \begin{smallmatrix}
  \vartheta_{2^{i-1}}(x) & 0_{2^{i-1}} \\
  0_{2^{i-1}} & 0_{2^{i-1}}
  \end{smallmatrix}
  }
  \\
  &=&
  \rukl{
  \begin{smallmatrix}
  \vartheta_{2^i}(x) & 0_{2^i} \\
  0_{2^i} & 0_{2^i}
  \end{smallmatrix}
  }
  = \rho_{2^i,2^{i+1}} \circ \vartheta_{2^i} (x).
\eglnoz
\eproof
We conclude that the isomorphisms $\vartheta_{2^i}$ give rise to an isomorphism between $\bar \cQ_2 = \ilim \gekl{\cQ_2; \Ad(s_2)}$ and $\ilim_i \gekl{M_{2^i}(\cQ_2); \rho_{2^i,2^{i+1}}}$. As the second inductive limit can be identified with $\cK(\ell^2(\Nz_0)) \otimes \cQ_2$, we arrive at the following result:
\bprop
\label{Q2KQ2}
There is a canonical isomorphism $\bar \cQ_2 \cong \cK(\ell^2(\Nz_0)) \otimes \cQ_2$.
\eprop

So the dilation $\bar \cQ_2$ is nothing else but the stabilization of $\cQ_2$. This result can also be obtained from classification theory (see \cite{Ror}, Chapter~8) because both $\cQ_2$ and $\bar \cQ_2$ are Kirchberg algebras satisfying the UCT by Remark~\ref{nuclear-UCT}. However, the point we make is that we get a canonical identification of the stabilization of $\cQ_2$ with $\bar \cQ_2$ which can be described explicitly.

Moreover, recall from Lemma~\ref{M'} that $(C_0(\Qz_2), \beta^{\operatorname{aff}}, \Zz[\halb] \rtimes \spkl{2})$ is the minimal automorphic dilation of $(C(\Zz_2), \alpha^{\operatorname{aff}}, \Zz\rtimes [2\rangle)$. By Proposition~\ref{Q2semigpcp}, the latter dynamical system (or rather semisystem) gives rise to a crossed product which can be identified with $\cQ_2$. So the crossed product obtained from the first dynamical system can be thought of as a dilation of $\cQ_2$ as well, just as $\bar \cQ_2$. In fact, the situation is precisely of the type considered in \cite{LarLi}.

\bprop
\label{barQ2Q2}
There is a canonical isomorphism $\bar \cQ_2 \cong C_0(\Qz_2) \rtimes_{\beta^{\operatorname{aff}}} (\Zz[\halb] \rtimes \spkl{2})$.
\eprop

\bproof
The assertion follows from \cite[Theorem~1.3]{LarLi} because $[2\rangle$ is cofinal in $\Zz\rtimes [2\rangle)$. To apply \cite[Theorem~1.3]{LarLi}, just proceed in the same way as in \cite[Section~2.1]{LarLi}.
\eproof

\bremark\label{nuclear-UCT}
Corollary~\ref{Q2MQ2} combined with Propositions~\ref{Q2KQ2} and \ref{barQ2Q2}
imply that $\cQ_2$ is stably isomorphic to
$C_0(\Qz_2)
\rtimes_{\beta^{\operatorname{aff}}} (\Zz[\halb] \rtimes \spkl{2})$. Since
$(\Zz[\halb] \rtimes \spkl{2})$ is amenable, the crossed product
$C_0(\Qz_2) \rtimes_{\beta^{\operatorname{aff}}} (\Zz[\halb] \rtimes \spkl{2})$
is nuclear and satisfies the UCT (see \cite{Ror}, Chapter~2). Hence  we see again that $\cQ_2$ is nuclear and
satisfies the UCT.
\eremark

\section{The 2-adic duality theorem}
\label{dualthm}

As explained in Section~\ref{croprodes}, we can think of $\cQ_2$ as a 2-adic version of the ring $C^*$-algebra of the integers. Given this interpretation, a natural question is whether the duality theorem from \cite[\S~3 and \S~4]{CuLi2} has a 2-adic analogue as well.

Our goal in this section is to answer this question by establishing a 2-adic duality result.

The semidirect product $\Zz[\halb] \rtimes \spkl{2}$ acts on $C_0(\Rz)$ and $C_0(\Qz_2)$ via affine transformations $\lambda^{\operatorname{aff}}$ and $\beta^{\operatorname{aff}}$, respectively, where both are given by
\bgl
\label{def-beta-generic}
  {(b,a)} \cdot f=f(a^{-1} (\sqcup-b))
\egl
for $f$ in $C_0(\Rz)$ or $C_0(\Qz_2)$, respectively. Our goal in this section is to prove that the corresponding crossed products $C_0(\Rz) \rtimes_{\lambda^{\operatorname{aff}}} (\Zz[\halb] \rtimes \spkl{2})$ and $C_0(\Qz_2) \rtimes_{\beta^{\operatorname{aff}}} (\Zz[\halb] \rtimes \spkl{2})$ are strongly Morita equivalent. The strategy of the proof is the same as in \cite[\S~3 and \S~4]{CuLi2}. We shall go through the details of the computations because we need to describe the imprimitivity bimodule which implements this equivalence.

We shall need to find the Pontryagin dual of $\Zz[\halb]$. To this end, we fix the model $\Zz_2 = \plim_i \gekl{\Zz / 2^i \Zz}$ as in the proof of Proposition \ref{DCZ2}:
\bgloz
  \Zz_2 \cong \menge{\zf =(z_i)_i \in \prod_{i=0}^{\infty} \Zz / 2^i \Zz}{z_i = p_{i+1,i}(z_{i+1})}
\egloz
where $p_{i+1,i}$ are the canonical projections $\Zz / 2^{i+1} \Zz \onto \Zz / 2^i \Zz$. The canonical projection $\Zz_2 \to \Zz / 2^i \Zz$ is $p_i$.

As $\Zz[\halb] / \Zz$ is the inductive limit of $2^{-i} \Zz /\Zz$, the Pontryagin dual of $\Zz[\halb] / \Zz$ is $\plim_i \gekl{\Zz / 2^i \Zz}$. The pairing $\spkl{n,t} =\e(nt)$ of $\Zz$ with $\Rz/\Zz$ identifies $\widehat{2^{-i}\Zz/\Zz}$ with $\Zz/2^i\Zz$, and carries the inclusion $2^{-i} \Zz / \Zz \into 2^{-(i+1)} \Zz /\Zz$ to the projection $ p_{i+1,i}$. Thus it implements an isomorphism
\bgl
\label{dual_of_Z2}
  \Zz_2 \to \widehat{\Zz[\halb] / \Zz}; \ \zf = (z_i) \ma \eckl{\tfrac{l}{2^i} + \Zz \ma \e(\tfrac{l}{2^i} z_i)}.
\egl

Let $p$ be the composition of the quotient map $\Qz_2 \onto \Qz_2/\Zz_2$ with the isomorphism $\Qz_2/\Zz_2 \cong \Zz[\halb] / \Zz$ determined by the inclusion $\Zz[\halb] \into \Qz_2$ (this inclusion gives rise to an isomorphism $\Zz[\halb] / \Zz \cong \Qz_2/\Zz_2$ because of $\Qz_2 = \Zz[\halb] + \Zz_2$ and $\Zz_2 \cap \Zz[\halb] = \Zz$). The map $t \ma \e(t)$ takes $\Zz[\halb] / \Zz$ into the subgroup $\Zz(2^\infty)$ of $\Tz$ consisting of all $\omega$ such that $\omega^{2^i}=1$ for $i \in \Nz$. In fact, $\widehat{\Zz[\halb] / \Zz}=\Zz(2^\infty)$, and the dual of $\Zz(2^\infty)$ is $\Zz_2$.
Then $\chi_0:\Qz_2 \to \Tz$, $\chi_0(\xf)=\e(p(\xf))$ is a character on $\Qz_2$ which implements the self-duality $\Qz_2\to \widehat{\Qz}_2$ sending $\xf \in \Qz_2$ to $\yf \ma \chi_0(\xf \yf)$. The restriction of this self-duality map to $\Zz_2$ is precisely the isomorphism from \eqref{dual_of_Z2}, i.e. we have $\e(\tfrac{l}{2^i} z_i) = \e(p(\tfrac{l}{2^i} \zf)) = \chi_0(\tfrac{l}{2^i} \zf)$ for every $\zf = (z_i)$ in $\Zz_2$ and every $\tfrac{l}{2^i} \in \Zz[\halb]$ (the crucial point is that $\tfrac{l}{2^i} z_i + \Zz = p(\tfrac{l}{2^i} \zf)$).

To describe the Pontryagin dual of $\Zz[\halb]$, consider the locally compact group $\Rz \times \Qz_2$ and its diagonal subgroup $\Delta \defeq \menge{(b,b)}{b \in \Zz[\halb]}$. Using $\Qz_2 = \Zz[\halb] + \Zz_2$ and $\Zz_2 \cap \Zz[\halb] = \Zz$, we obtain as above that the inclusion $\Rz \times \Zz_2 \into \Rz \times \Qz_2$ induces an isomorphism of topological groups
\bgloz
  \Rz \times \Zz_2 / \Delta \cap (\Rz \times \Zz_2) \cong \Rz \times \Qz_2 / \Delta.
\egloz
Note that $\Delta \cap (\Rz \times \Zz_2) = \menge{(b,b)}{b \in \Zz}$. The \emph{2-adic solenoid}
\bgloz
  \Sigma_2 \defeq \Rz \times \Zz_2 / \Delta \cap (\Rz \times \Zz_2) \cong \Rz \times \Qz_2 / \Delta
\egloz
is a compact connected abelian group, see for example \cite[Theorem~10.13]{HeRo} for the general case of the full solenoid. We let $\pi: \Rz \times \Qz_2 \to \Rz\times \Qz_2 / \Delta$ be the quotient map. We shall often write $[r,\xf]=\pi(r,\xf)$ for $r\in \Rz$ and $\xf \in \Qz_2$. The next result can be proven either by following the
proof of \cite[Lemma~3.5]{CuLi2} using the description of $\widehat{\Zz[\halb] / \Zz}$ developed above, or by specialising the arguments in \cite[\S~(25.3) and \S~(25.4)]{HeRo} to the 2-adic solenoid.

\blemma
The map $\chi:\Sigma_2 \rightarrow \widehat{\Zz[\halb]}$ given by
\bgl
\label{def_of_chi}
  \chi([r,\xf])= [b \ma \e(rb) \e(-p(\xf b))]
\egl
for $(r,\xf)\in \Rz\times \Qz_2$ implements a topological isomorphism of $\Sigma_2$ onto the Pontryagin dual of $\Zz[\halb]$.
\elemma

There is a subtle point regarding the choice of the character which implements the Pontryagin self-duality of $\Rz$. Since our actions involve simultaneously the self-dual groups $\Rz$ and $\Qz_2$, and since the self-duality of $\Qz_2$ is given by $\chi_0$, we shall choose the self-duality $\Rz \cong \widehat{\Rz}$ given by $t \ma \eckl{x \ma \e(tx)}$ (in other words, we allow a factor of $2 \pi$ in the exponential function).

Fix the embedding $\Zz[\halb] \into \widehat{\Rz}$ given by $b \ma \eckl{t \ma \e{(-tb)}}$. There are actions $\tau$ of $\Zz[\halb]$ on $C^*(\Rz)$ and $\sigma$ of $\Rz$ on $C^*(\Zz[\halb])$ given, respectively, by
\bgloz
  (\tau_b f)(t)=f(t) \e(-tb) \,\,\, \text{ and }\,\,\,(\sigma_t g)(b) = g(b) \e(bt),
\egloz
for $f\in C_c(\Rz)\subseteq C^*(\Rz)$, $g \in C_c(\Zz[\halb])\subseteq C^*(\Zz[\halb])$, $b \in \Zz[\halb]$ and $t \in \Rz$. Then the duality result \cite[Lemma~4.3]{CuLi2} applied to our situation (with $G = \Rz$ and $H = \Zz[\halb]$) gives an isomorphism
\bgloz
  \varphi: C^*(\Rz) \rtimes_{\tau} H \rightarrow C^*(H) \rtimes_{\sigma} \Rz
\egloz
which sends an element $f\in C_c(H, C_c(\Rz))$ to $\varphi(f)$ given by
\bgloz
  ((\varphi f)(t))(b)=\e(tb)(f(b))(t)
\egloz
for $t\in \Rz$ and $b \in \Zz[\halb]$.

The Fourier transform $f \ma \hat{f}$ where $\hat{f}(x)=\int_{\Rz} \e(tx) f(t) dt$ for $f\in C_c(\Rz)\subseteq C^*(\Rz)$ implements an isomorphism $F_\Rz: C^*(\Rz) \to C_0(\widehat{\Rz})$. Note that we identify $\Rz$ with $\widehat{\Rz}$ via $t \ma \eckl{x \ma \e(tx)}$ as explained before. We shall let $\check{f}$ denote the inverse Fourier transform of a function $f \in C_c(\Rz)\subseteq C_0(\Rz)$. Note that $F_\Rz$ carries the action $\tau$ of $\Zz[\halb]$ on $C^*(\Rz)$ into the action $\hat{\tau}$ on $C_0(\Rz)$ given by $(\hat{\tau}_b(g))(x)=g(x-b)$ for $g \in C_c(\Rz)$. Thus $F_\Rz$ implements
an isomorphism
\bgloz
  F_{\Rz}: C^*(\Rz) \rtimes_{\tau} \Zz[\halb] \to C_0(\Rz) \rtimes_{\hat{\tau}} \Zz[\halb].
\egloz

The Fourier transform on $\Zz[\halb]$ is $F_{\Zz[\halb]}:C^*(\Zz[\halb]) \to C(\widehat{\Zz[\halb]})$ given by
\bgloz
  (F_{\Zz[\halb]}g)(\gamma)=\sum_{b \in \Zz[\halb]} g(b) \gamma(b) \text{ for } g\in C_c(\Zz[\halb]) \subseteq C^*(\Zz[\halb]) \text{ and } \gamma \in \widehat{\Zz[\halb]}.
\egloz
Let $\chi^*: C(\widehat{\Zz[\halb]})\to C(\Sigma_2)$ be the isomorphism coming from the identification \eqref{def_of_chi}. Since
\bglnoz
  (F_{\Zz[\halb]}(\sigma_t(g)))(\chi([r,\xf]))
  &=& \sum_{b \in \Zz[\halb]} g(b) \e(tb) \e(rb) \e(-p(\xf b)) \\
  &=& \sum_{b \in \Zz[\halb]} g(b) \chi([t+r,\xf])(b),
\eglnoz
it follows that $\chi^*\circ F_{\Zz[\halb]}$ carries $\sigma$ into the action $\hat{\sigma}$ on $C(\Sigma_2)$ given by
\bgloz
  \hat{\sigma}_t (h)([r,\xf])= h([r+t,\xf])
\egloz
for $t \in \Rz$, $h \in C(\Sigma_2)$ and $[r,\xf] \in \Sigma_2$. Hence it gives rise to an isomorphism
\bgloz
  F_{\Zz[\halb],\chi}: C^*(\Zz[\halb]) \rtimes_{\sigma} \Rz \to C(\Sigma_2) \rtimes_{\hat{\sigma}} \Rz
\egloz
which takes $g\in C_c(\Rz, C_c(\Zz[\halb]))$ to the function in $C_c(\Rz,C(\Sigma_2))$ given by
\bgloz
  ((F_{\Zz[\halb],\chi}g)(t))([r,\xf]) = \sum_{d \in \Zz[\halb]} \chi([r,\xf])(d)g(t)(d).
\egloz

Now let $\varphi^F$ be the composition $F_{\Zz[\halb],\chi} \circ \varphi \circ F_\Rz^{-1}$. Then for $f\in C_c(\Zz[\halb], C_c(\Rz))$ we have
\bgln
  ((\varphi^F f)(t))([r, \xf])
  &=& \sum_{d \in \Zz[\halb]} \chi([r,\xf])(d)(\varphi \circ F_\Rz^{-1}(f))(t)(d) \nonumber \\
  &=& \sum_{d \in \Zz[\halb]} \chi([r,\xf])(d) \e(td) ((F_\Rz^{-1} f)(d))(t) \nonumber \\
  &=& \sum_{d \in \Zz[\halb]} \e((r+t)d) \e(-p(\xf d)) (F_\Rz^{-1} f)(d))(t). \label{formula-phi^F}
\egln

Next we claim that $\varphi^F$ is equivariant for the actions of $\spkl{2}$ on $C_0(\Rz) \rtimes_{\hat{\tau}} \Zz[\halb]$ and on $C(\Sigma_2) \rtimes_{\hat{\sigma}} \Rz$ given, respectively, by
\bgloz
(\rho_a (f)(b))(t)=f(a^{-1} b)(a^{-1} t)
\egloz
and
\bgloz
(\ti{\rho}_a (h))(t)([r,\xf]) = a h(ta)([ra,\xf a])
\egloz
for $a \in \spkl{2}$, $f\in C_c(\Zz[\halb], C_c(\Rz))$, $h\in C_c(\Rz, C(\Sigma_2))$, $b \in \Zz[\halb]$, $t \in \Rz$ and $[r,\xf] \in \Sigma_2$. Indeed,
\bglnoz
  && (\varphi^F (\rho_a (f))(t))([r,\xf]) \\
  &=& \sum_{d \in \Zz[\halb]} \e((t+r)d) \e(-p(\xf d)) \int_{\Rz} \e(-ts) f(a^{-1} d)(a^{-1} s) ds \\
  &=& a \sum_{d \in \Zz[\halb]} \e((t+r)ad) \e(-p(\xf ad)) (F_{\Rz}^{-1} f (d))(ta) \\
  &=& a ((\varphi^F f)(ta))([ra,\xf a]) \text{ by } \eqref{formula-phi^F} \\
  &=& (\ti{\rho}_a (\varphi^F f))(t)([r,\xf]).
\eglnoz

Hence $\varphi^F$ induces an isomorphism
$(C_0(\Rz)\rtimes_{\hat{\tau}} \Zz[\halb]) \rtimes_{\rho} \spkl{2} \cong (C(\Sigma_2) \rtimes_{\hat{\sigma}} \Rz)\rtimes_{\ti{\rho}} \spkl{2}$.
We rewrite this by twice applying \cite[Proposition~3.11]{Wil}. The semidirect product $\Rz \rtimes \spkl{2}$ arising from the action $a \cdot t = t a^{-1}$ for $a \in \spkl{2}$ and $t \in \Rz$ acts on $C(\Sigma_2)$ by
\bgl
\label{action-of-R-and-2}
  (\hat{\sigma} \rtimes \ti{\rho})_{(t,a)} (h)([r,\xf]) = h([(t+r)a, \xf a])
\egl
for all $(t,a) \in \Rz \rtimes \spkl{2}$ and $h \in C(\Sigma_2)$. Since
\bgloz
  (\hat{\sigma} \rtimes \ti{\rho})_{(0,a)} (h(ta)) ([r,\xf]) = h(ta)([ra,\xf a])= a^{-1} (\ti{\rho}_a (h))(t)([r,\xf])
\egloz
for $h \in C_c(\Rz, C(\Sigma_2))$, \cite[Proposition~3.11]{Wil} gives an isomorphism
\bgloz
  (C(\Sigma_2) \rtimes_{\hat{\sigma}} \Rz) \rtimes_{\ti{\rho}} \spkl{2}
  \cong C(\Sigma_2) \rtimes_{\hat{\sigma} \rtimes \ti{\rho}} (\Rz \rtimes \spkl{2}).
\egloz
This isomorphism sends $h$ in $C_c(\spkl{2},C_c(\Rz,C(\Sigma_2))) \subseteq (C(\Sigma_2) \rtimes_{\hat{\sigma}} \Rz) \rtimes_{\ti{\rho}} \spkl{2}$ to the function
\bgloz
  \Rz \rtimes \spkl{2} \ni (t,a) \ma a^{-1} h(a)(t) \in C(\Sigma_2)
\egloz
in $C_c(\Rz \rtimes \spkl{2}, C(\Sigma_2)) \subseteq C(\Sigma_2) \rtimes_{(\hat{\sigma} \rtimes \ti{\rho})} (\Rz \rtimes \spkl{2})$.

The group $\Zz[\halb]\rtimes \spkl{2}$ acts on $C_0(\Rz)$ by affine transformations $\lambda^{\operatorname{aff}}$ given in \eqref{def-beta-generic}.
But
\bgloz
  \lambda^{\operatorname{aff}}_{(0,a)}(g(a^{-1}b))(x)=(g(a^{-1}b))(a^{-1}x)=(\rho_a(g)(b))(x)
\egloz
for all $g \in C_c(\Zz[\halb], C_c(\Rz))$, hence \cite[Proposition~3.11]{Wil} yields
\bgloz
  (C_0(\Rz)\rtimes_{\hat{\tau}} \Zz[\halb]) \rtimes_{\rho} \spkl{2} \cong C_0(\Rz) \rtimes_{\lambda^{\operatorname{aff}}} (\Zz[\halb] \rtimes \spkl{2}).
\egloz
The isomorphism sends a function $g$ in $C_c(\spkl{2},C_c(\Zz[\halb],C_0(\Rz)))$ to the function
\bgloz
  \Zz[\halb] \rtimes \spkl{2} \ni (b,a) \ma g(a)(b) \in C_0(\Rz)
\egloz
in $C_c(\Rz \rtimes \spkl{2}, C(\Sigma_2)) \subseteq C_0(\Rz) \rtimes_{\lambda^{\operatorname{aff}}} (\Rz \rtimes \spkl{2})$.

Putting all this together, we have the following
\bprop
\label{d'}
There is an isomorphism
\bgloz
  \Phi: C_0(\Rz) \rtimes_{\lambda^{\operatorname{aff}}} (\Zz[\halb]\rtimes \spkl{2}) \to C(\Sigma_2) \rtimes_{\hat{\sigma} \rtimes \ti{\rho}} (\Rz\rtimes \spkl{2})
\egloz
which sends $f \in C_c(\Zz[\halb]\rtimes \spkl{2}, C_c(\Rz))$ to the function in $C_c(\Rz \rtimes \spkl{2}, C(\Sigma_2))$ given by
\bgl
\label{Phi-iso}
  \Phi(f)(t,a)([r,\xf])= a^{-1} \sum_{b \in \Zz[\halb]} \e((r+t)b) \e(-p(\xf b)) F_{\Rz}^{-1}(f(b,a))(t).
\egl
for $(t,a) \in \Rz \rtimes \spkl{2}$ and $[r,\xf]\in \Sigma_2$.
\eprop

The next step towards the 2-adic duality theorem is to note that the action $\hat{\sigma} \rtimes \ti{\rho}$ from \eqref{action-of-R-and-2} is implemented by the action of $\Rz \rtimes \spkl{2}$ on $\Sigma_2$ given by
\bgl
  \label{first-transformation-gpd}
  [r,\xf] \cdot (t,a) \defeq [(r+t)a,\xf a]
\egl
for all $(t,a) \in \Rz \rtimes \spkl{2}$ and $[r,\xf] \in \Sigma_2$. Before we present the next result, we specify our convention for the source and range maps on a transformation groupoid $X \rtimes G$, namely we let the source of $(x,g)$ be $x \cdot g$ and the range of $(x,g)$ be $x$.

\bprop
\label{d''}
The transformation groupoids $\cG$ and $\ti{\cG}$ associated, respectively, with the right actions of $\Rz \rtimes \spkl{2}$ on $\Rz \times \Qz_2 / \Delta$ from
\eqref{first-transformation-gpd} and of $\Zz[\halb] \rtimes \spkl{2}$ on $\Qz_2$ via
\bgloz
  \xf \cdot (b,a) = (\xf-b)a,
\egloz
are equivalent in the sense of Muhly-Renault-Williams \cite{MRW}.
\eprop

\bproof
The proof is analogous to the proof of Proposition~3.8 in \cite{CuLi2} and of Lemma~4.4 in \cite{CuLi2}.

Let $N \defeq \pi(\{0\}\times \Zz_2) \subseteq \cG^0$ and $\ti{N} \defeq \Zz_2 \subseteq \ti{\cG}^0$. Both $N$ and $\ti{N}$ fulfill the conditions of \cite[Example~2.7]{MRW}: For example, $\Zz_2$ is closed (in fact compact) and open in $\Qz_2$, $\ti{\cG}^0 = \bigcup_{(b,a) \in \Zz[\halb]\rtimes \spkl{2}} \Zz_2 \cdot (b,a)$ and the restricted range and source maps $\ti{r} \vert_{\ti{\cG}_{\ti{N}}}: \ti{\cG}_{\ti{N}} \ri \ti{\cG}^0$, $\ti{s} \vert_{\ti{\cG}_{\ti{N}}}^{\ti{N}}: \ti{\cG}_{\ti{N}} \ri \ti{N}$ are open, so
\bgloz
  \ti{\cG}_{\ti{N}} = \menge{(\xf,(b,a)) \in \ti{\cG}}{(\xf -b)a \in \Zz_2}
\egloz
implements an equivalence between $\ti{\cG}$ and $\ti{\cG}_{\ti{N}}^{\ti{N}}$ given by
\bgloz
  \ti{\cG}_{\ti{N}}^{\ti{N}} = \menge{(\zf,(b,a)) \in \ti{\cG}}{\zf \in \Zz_2 \text{ and } (\zf -b)a \in \Zz_2}.
\egloz
Similarly, $N$ is a closed subset of $\cG^0 = \Sigma_2$, $\cG^0 = \bigcup_{(t,a) \in \Rz\rtimes \spkl{2}} N \cdot (t,a)$ and the restricted range and source maps $r \vert_{\cG_N}: \cG_N \ri \cG^0$, $s \vert_{\cG_N}^{N}: \cG_N \ri N$ are open, so
\bgloz
  \cG_N = \menge{([r,\xf],(t,a)) \in \cG}{[(r+t)a,\xf a] \in N}
\egloz
implements an equivalence between $\cG$ and $\cG_N^N$ given by
\bgloz
  \cG_N^N = \menge{([0, \zf],(t,a)) \in \cG}{\zf \in \Zz_2 \text{ and } [ta,\zf a] \in N}.
\egloz

 Note here that $([0,\zf],(t,a)) \in \cG_N^N$ entails that $[ta, \zf a] =[0,\zf']$ for some $\zf' \in \Zz_2$, which therefore implies $ta = \zf a- \zf' \in \Zz[\halb]$, and then necessarily $t \in \Zz[\halb]$. In addition, we see that $[ta, \zf a]=[0, \zf a-ta]$.

It follows that the map $(\zf,(b,a)) \ma ([0,\zf],(b,a))$ is an isomorphism from $\ti{\cG}_{\ti{N}}^{\ti{N}}$ onto $\cG_N^N$, seen as topological groupoids. Thus $ \cG \sim \cG_N^N \cong \ti{\cG}_{\ti{N}}^{\ti{N}} \sim \ti{\cG}$, as claimed.
\eproof

We record for later use that the isomorphism of $\ti{\cG}_{\ti{N}}^{\ti{N}}$ onto $\cG_N^N$ gives rise to an isomorphism $\Psi_{N, \ti{N}}$ from $C^*(\cG_N^N)$ onto $C^*(\ti{\cG}_{\ti{N}}^{\ti{N}})$.

\bcor
\label{def-Psi}
There is an isomorphism $\Psi:C^*(\ti{\cG}_{\ti{N}}^{\ti{N}})\to \cQ_2$ which maps $\1z_{\Zz_2\times \{(1,1)\}}\in C_c(\ti{\cG}_{\ti{N}}^{\ti{N}})$ to $u$ and $\1z_{2 \Zz_2 \times \{(0, 2^{-1})\}} \in C_c(\ti{\cG}_{\ti{N}}^{\ti{N}})$ to $s_2$.
\ecor

\bproof
Note that $C^*(\ti{\cG}_{\ti{N}}^{\ti{N}})$ is isomorphic to the corner $\1z_{\Zz_2} C^*(\ti{\cG}) \1z_{\Zz_2}$ of $C^*(\cG)$ given by the full projection $\1z_{\Zz_2}$. The algebra $C^*(\ti{\cG})$ is isomorphic to $C_0(\Qz_2) \rtimes_{\ti{\beta}} (\Zz[\halb] \rtimes \spkl{2})$ arising from the action $\ti{\beta}$ of $\Zz[\halb]\rtimes\spkl{2}$ given by $(b,a) \cdot f=f((\sqcup-b)a)$ for $(b,a) \in \Zz[\halb] \rtimes \spkl{2}$ and $f \in C_c(\Qz_2)$. This action is not
$\beta^{\operatorname{aff}}$, however, under the isomorphism $\spkl{2} \ni a \ma a^{-1} \in \spkl{2}$
we have
\bgloz
C_0(\Qz_2) \rtimes_{\ti{\beta}} (\Zz[\halb] \rtimes \spkl{2})\cong C_0(\Qz_2) \rtimes_{\beta^{\operatorname{aff}}} (\Zz[\halb] \rtimes \spkl{2}).
\egloz
Thus by Lemma~\ref{M'}, $C^*(\ti{\cG}_{\ti{N}}^{\ti{N}})$ is isomorphic to $C(\Zz_2)\rte_{\alpha^{\operatorname{aff}}} (\Zz\rtimes [2\rangle)$, where the isomorphism takes $\1z_{\Zz_2\times \{(1,1)\}}$ to $w_{(1,1)}$ and $\1z_{2\Zz_2\times\{(0, 2^{-1})\}}$ to $w_{(0,2)}$. Applying Corollary~\ref{Q2Z2} finishes the proof.
\eproof

The following is the second main result of this paper:

\btheo[The 2-adic duality theorem]
\label{d}
Let $\beta^{\operatorname{aff}}$ and $\lambda^{\operatorname{aff}}$ be the actions of $\Zz[\halb] \rtimes \spkl{2}$ by affine transformations on $C_0(\Qz_2)$ and $C_0(\Rz)$, respectively, given by \eqref{aff-transf-Q2} and \eqref{action-of-R-and-2}. Then $C_0(\Rz) \rtimes_{\lambda^{\operatorname{aff}}} (\Zz[\halb] \rtimes \spkl{2})$ and $C_0(\Qz_2) \rtimes_{\beta^{\operatorname{aff}}} (\Zz[\halb] \rtimes \spkl{2})$ are Morita equivalent.
\etheo

\bproof
Proposition \ref{d''} and \cite[Theorem~2.8]{MRW} yield that $C^*(\cG)$ and $C^*(\ti{\cG})$ are Morita equivalent.

The first groupoid $C^*$-algebra can be identified with $C(\Sigma_2) \rtimes_{\hat{\sigma} \rtimes \ti{\rho}} (\Rz \rtimes \spkl{2})$ from Proposition \ref{d'}. By the proof of Corollary~\ref{def-Psi}, the second groupoid $C^*$-algebra can be identified with $C_0(\Qz_2) \rtimes_{\ti{\beta}} (\Zz[\halb] \rtimes \spkl{2})$, and thus with $C_0(\Qz_2) \rtimes_{\beta^{\operatorname{aff}}} (\Zz[\halb] \rtimes \spkl{2})$ using the isomorphism $\spkl{2} \ni a \ma a^{-1} \in \spkl{2}$ as in the proof of Corollary~\ref{def-Psi}. Hence we have
\begin{align*}
  C_0(\Rz) \rtimes_{\lambda^{\operatorname{aff}}} (\Zz[\halb] \rtimes \spkl{2})
  &\overset{\Phi}{\cong} C(\Sigma_2) \rtimes_{\hat{\sigma} \rtimes \ti{\rho}} (\Rz \rtimes \spkl{2}) \\
  &\sim_M C_0(\Qz_2) \rtimes_{\beta^{\operatorname{aff}}} (\Zz[\halb] \rtimes \spkl{2}).
\end{align*}
This proves the 2-adic duality theorem.
\eproof

From now on, we write $C_0(\Rz) \rtimes \Zz[\halb] \rtimes \spkl{2}$ and $C_0(\Qz_2) \rtimes \Zz[\halb] \rtimes \spkl{2}$ for the crossed products arising from the actions $\beta^{\operatorname{aff}}$ and $\lambda^{\operatorname{aff}}$, respectively.

\bremark\label{Zhang-dichotomy}
The 2-adic duality theorem implies that $C_0(\Rz) \rtimes \Zz[\halb] \rtimes \spkl{2}$ and $C_0(\Qz_2) \rtimes \Zz[\halb] \rtimes \spkl{2}$ are actually isomorphic, not just Morita equivalent, as we shall now explain.

It follows from Theorem \ref{d} that these crossed products are stably isomorphic because both of them are separable. Moreover, $\cQ_2$ is purely infinite and simple (see Theorem \ref{pis}), thus $C_0(\Qz_2) \rtimes \Zz[\halb] \rtimes \spkl{2}$ has this property as well (see Corollary \ref{Q2MQ2} and \cite{Ror}, Proposition~4.1.8~(i)). But then, $C_0(\Rz) \rtimes \Zz[\halb] \rtimes \spkl{2}$ must also be purely infinite and simple (see \cite{Ror}, Proposition~4.1.8~(i); $C_0(\Rz) \rtimes \Zz[\halb] \rtimes \spkl{2}$ and $C_0(\Qz_2) \rtimes \Zz[\halb] \rtimes \spkl{2}$ are stably isomorphic as explained above). Now we can apply \an{Zhang's dichotomy} (see \cite{Ror}, Proposition~4.1.3) which yields that both $C_0(\Rz) \rtimes \Zz[\halb] \rtimes \spkl{2}$ and $C_0(\Qz_2) \rtimes \Zz[\halb] \rtimes \spkl{2}$ must be stable as they are not unital. So we have seen that these two crossed products are stably isomorphic and that they are stable, hence they must be isomorphic.

But of course, these arguments do not yield a concrete isomorphism between
$C_0(\Rz) \rtimes \Zz[\halb] \rtimes \spkl{2}$ and $C_0(\Qz_2) \rtimes \Zz[\halb] \rtimes \spkl{2}$.
\eremark

\section{A concrete bimodule}

Corollary \ref{Q2MQ2} showed that $\cQ_2$ and $C_0(\Qz_2) \rtimes \Zz[\halb] \rtimes \spkl{2}$ are Morita equivalent. Thus by Theorem~\ref{d}, $\cQ_2$ is Morita equivalent to $C_0(\Rz) \rtimes \Zz[\halb] \rtimes \spkl{2}$. There are two canonical representations around: The first is the (restricted) left regular representation $\lambda_2$ of $\cQ_2$ on $\ell^2(\Zz)$ constructed in Section \ref{Q2}. Moreover, the unitary representation of $\Zz$ on $\Ell^2(\Rz)$ given by $(n \cdot \xi)(x)= \xi(x-n)$ for $n \in \Zz$, $x \in \Rz$ extends to a unitary representation $T: \Zz[\halb] \to \cU(\Ell^2(\Rz))$ given by $(T_b \xi)(x)= \xi(x-b)$ for $b \in \Zz[\halb]$. Let $D$ be the unitary on $\Ell^2(\Rz)$ given by $(D \xi)(x)=2^{-1/2} \xi(2^{-1}x)$. This unitary gives rise to an action of $\Zz$ on $\Ell^2(\Rz)$ which we write as an action of $\spkl{2}$ in the following way: Given $a=2^i$ with $i \in \Zz$, put $D_a \defeq D^i$ (if $i < 0$, then of course $D^i$ means $(D^{-1})^{\abs{i}} = (D^*)^{\abs{i}}$). Then
$D_a T_b = T_{ab} D_a$ for all $b \in \Zz[\halb]$ and $a \in \spkl{2}$, and so
\bgloz
  (T\rtimes D)_{(b,a)} = T_b D_a \text{ for } (b,a)\in \Zz[\halb] \rtimes \spkl{2}
\egloz
is a unitary representation of $\Zz[\halb] \rtimes \spkl{2}$ on $\Ell^2(\Rz)$. This is a classical
  example of a representation of $\Zz[\halb] \rtimes \spkl{2}$ (which is known as a Baumslag-Solitar group), see e.g. \cite{Dut}.
  Let $M$ be the representation of $C_0(\Rz)$ on $\Ell^2(\Rz)$ via multiplication operators. Then $(M, T \rtimes D)$ is a covariant pair for the dynamical system $(C_0(\Rz), \Zz[\halb] \rtimes \spkl{2}, {\lambda^{\operatorname{aff}}})$ and thus gives rise to a representation of $C_0(\Rz) \rtimes \Zz[\halb] \rtimes \spkl{2}$ on $\Ell^2(\Rz)$. The second representation of interest here is
\bgl
  \label{standard-rep}
  \pi(f u_{(b,a)}) \defeq M_f (T\rtimes D)_{(b,a)},
\egl
for all $f \in C_0(\Rz)$, $(b,a) \in \Zz[\halb] \rtimes \spkl{2}$ and $\xi \in \Ell^2(\Rz)$. In other words,
\bgloz
  (\pi(f \otimes \1z_{(b,a)}) \xi)(x)=f(x)a^{-1/2}\xi(a^{-1}(x-b))
\egloz
for $f \in C_0(\Rz)$, $(b,a)\in \Zz[\halb] \rtimes \spkl{2}$ and $\xi \in \Ell^2(\Rz)$, where we view $f \otimes \1z_{(b,a)} \in C_c(\Zz[\halb] \rtimes \spkl{2},C_0(\Rz))$ as an element in $C_0(\Rz) \rtimes \Zz[\halb] \rtimes \spkl{2}$.

Since Morita equivalence yields a $1$--$1$ correspondence between representations of the $C^*$-algebras, a natural question is whether the representations $\lambda_2$ and $\pi$ correspond to one another with respect to the imprimitivity bimodule which we obtain from the proof of the 2-adic duality theorem.

To answer this question, we take the $C_0(\Rz)\rtimes \Zz[\halb]\rtimes\spkl{2}$--$\cQ_2$-imprimitivity bimodule $X$ given by the proof of the 2-adic duality theorem and use it to induce representations. We obtain $X$ as a completion of a left-$R$, right-$Q$ pre-imprimitivity bimodule $X_0$ for appropriate dense subalgebras $R$ of $C_0(\Rz)\rtimes \Zz[\halb]\rtimes\spkl{2}$ and $Q$ of $\cQ_2$, respectively.

We can describe our situation by a diagram, as follows:
$$
  \xymatrix{C^*(\cG)\ar[d]_{\Phi_{\text{gpd},\rtimes}} \ar@{--}[r] &{Z} \ar@{--}[r] & C^*({\cG}_N^N) \ar[d]^{\Psi_{N, \ti{N}}}\\
  C(\Sigma_2)\rtimes(\Rz\rtimes\spkl{2}) \ar[d]_{\Phi^{-1}} & {} & C^*(\ti{\cG}^{\ti{N}}_{\ti{N}})\ar[d]^{\Psi}\\
  C_0(\Rz)\rtimes(\Zz[\halb]\rtimes\spkl{2}) \ar@{--}[r]&{X} \ar@{--}[r]&\cQ_2}
$$

Here, the imprimitivity bimodule $Z$ is obtained from the bimodule which implements the equivalence of $C^*(\cG)$ and $C^*({\cG}_N^N)$ as explained in \cite{MRW}. By composing with $\Phi^{-1} \circ \Phi_{\text{gpd},\rtimes}$ we alter the algebra acting from the left and by composing with $\Psi\circ \Psi_{N, \ti{N}}$ we alter the
algebra acting from the right. Here $\Phi_{\text{gpd},\rtimes}: C^*(\cG) \cong C(\Sigma_2) \rtimes_{\hat{\sigma} \rtimes \ti{\rho}} (\Rz \rtimes \spkl{2})$ is the canonical isomorphism identifying the groupoid $C^*$-algebra of a transformation groupoid with the corresponding crossed product. The result will be a $R$--$Q$ pre-imprimitivity bimodule $X_0$, where $R:=\Phi^{-1}(C_c(\Sigma_2, \Rz\rtimes\spkl{2}))$ is a dense subalgebra of $C_0(\Rz)\rtimes(\Zz[\halb]\rtimes\spkl{2})$ and $Q$ is the dense subalgebra of $\cQ_2$ generated by $s_2$ and $u$.

Let us start with the underlying vector space of $X_0$. The left-$C_c(\cG)$ and right-$C_c({\cG}_N^N)$ pre-imprimitivity bimodule which implements the equivalence of $C^*(\cG)$ and $C^*({\cG}_N^N)$ is based on the space $C_c(\cG_N)$. In order to describe the actions and the inner-products, note first that the range map of $\cG$ is
\bgloz
  r: \cG \to \Sigma_2; ([r,\zf],(t,a)) \ma [r,\zf]
\egloz
and the source map is
\bgloz
  s: \cG \to \Sigma_2; ([r,\zf],(t,a)) \ma [(r+t)a,\zf a].
\egloz
The range and source maps of $\cG_N^N$ are obtained as the restrictions of the range and source maps of $\cG$. On $\cG$ we take the Haar system $\lambda^{[r,\xf]} = \delta_{[r, \xf]} \times \lambda$ where $\lambda$ is a left Haar measure on $\Rz \rtimes \spkl{2}$, i.e. $\int_{\Rz \rtimes \spkl{2}} f d \lambda = \sum_{a \in \spkl{2}} \int_\Rz f(t,a) a dt$.
On $\cG_N^N$ we take the Haar system given by $\mu^{[0,\zf]} = \delta_{[0,\zf]} \times \mu$ where $\mu$ is the canonical Haar measure on $\Zz[\halb] \rtimes \spkl{2}$ (so $\mu$ is simply given by the counting measure).

Let $\cG^N \defeq \menge{([r,\xf],(t,a)) \in \cG}{[r,\xf] \in N = \pi(\gekl{0} \times \Zz_2)}$. Then the inverse map $\cG^N \to \cG_N$, $\gamma \ma \gamma^{-1}$ is an isomorphism. The isomorphism $\Zz_2 \ni \zf \ma [0,\zf] \in N$ gives rise to an isomorphism from $\Zz_2 \times (\Rz \rtimes \spkl{2})$ onto $\cG^N$, and following this with the inverse map induces an isomorphism $\Theta: C_c(\cG_N) \to C_c(\Zz_2 \times (\Rz \rtimes \spkl{2}))$ such that
\bgl
\label{def-theta}
  (\Theta \phi)(\zf,(t,a)) = \phi(([0,\zf],(t,a))^{-1}) = \phi([ta,\zf a],(-ta,a^{-1}))
\egl
for $\phi \in C_c(\cG_N)$. Let
\bgloz
X_0 \defeq \Theta(C_c(\cG_N)) = C_c(\Zz_2 \times (\Rz \rtimes \spkl{2})),
\egloz
and note that
\bgl
\label{dense-lin-span}
  \lspan \menge{\1z_{l+k\Zz_2} \otimes \xi \otimes \1z_{\gekl{m}}}{k \in [2\rangle, l \in \Zz, \xi \in C_c(\Rz), m \in \spkl{2}} \text{ is dense in } X_0
\egl
with respect to the inductive limit topology.

Now let us determine the left action of $R$ on $X_0$. First, $C_c(\cG)$ acts on the left on $X_0$ by $g\cdot \Theta(\phi)=\Theta(g{\bullet} \phi)$ for $g\in C_c(\cG)$
and $\phi \in C_c(\cG_N)$, where $\bullet$ refers to right and left actions described in \cite[Page 11]{MRW}. The formula for the left action is
\bgl
\label{MRW-leftaction}
  (g \bullet \phi)(\zeta) = \int_{\cG} g(\gamma) \phi(\gamma^{-1} \zeta) d \lambda^{r(\zeta)}(\gamma) \fuer \zeta \in \cG_N.
\egl
Ultimately, we want a left action of $R = \Phi^{-1}(C_c(\cG))$ and an $R$-valued inner product on $X_0$. It certainly suffices to determine the action for elements of the form $f \otimes \1z_{\gekl{(d,c)}} \in R$ with $f \in F_{\Rz}(C_c(\Rz))$ and $(d,c) \in \Zz[\halb] \rtimes \spkl{2}$ because the linear span of these elements is dense in $R$. As explained above, the left action of $R$ on $X_0$ is given by the left action of $C_c(\cG)$ on $X_0$ upon applying the isomorphism $\Phi_{\text{gpd},\rtimes}^{-1} \circ \Phi$ to the element that acts, i.e. we set
$(f \otimes \1z_{\gekl{(d,c)}}) \cdot \Theta(\phi) \defeq (\Phi_{\text{gpd},\rtimes}^{-1} \circ \Phi)(f \otimes \1z_{\gekl{(d,c)}}) \cdot \Theta(\phi)$ for $\phi \in C_c(\cG_N)$.

By \eqref{Phi-iso}, the function $\Phi(f \otimes \1z_{\gekl{(d,c)}})$ has support included in $\gekl{(\cdot,c)} \subseteq \Rz \rtimes \spkl{2}$, and
\bgloz
  \Phi(f \otimes \1z_{\gekl{(d,c)}})(t,a)([r,\xf])=c^{-1} \e((r+t)d) \e(-p(\xf d))\check{f}(t) \1z_{\gekl{c}}(a)
\egloz
and hence
\bgloz
  (\Phi_{\text{gpd},\rtimes}^{-1} \circ \Phi)(f \otimes \1z_{\gekl{(d,c)}})
([r,\xf], (t,a))
=c^{-\halb} \e((r+t)d) \e(-p(\xf d))\check{f}(t) \1z_{\gekl{c}}(a)
\egloz
since $\Phi_{\text{gpd},\rtimes}$ is the canonical identification $C^*(\cG) \cong C(\Sigma_2) \rtimes_{\hat{\sigma} \rtimes \ti{\rho}} (\Rz \rtimes \spkl{2})$ which sends $h$ in $C_c(\cG) = C_c(\Sigma_2 \times (\Rz \rtimes \spkl{2}))$ to the function
\bgloz
  \Rz \rtimes \spkl{2} \ni (t,a) \ma a^{-\halb} h(\sqcup,(t,a)) \in C(\Sigma_2).
\egloz
The factor $a^{-\halb}$ comes in because of the modular function associated with the group $\Rz \rtimes \spkl{2}$.

We therefore have
\bgln
  && ((f\otimes\1z_{\{(d,c)\}} )\cdot \Theta(\phi))(\zf,(t,a)) \nonumber \\
  &=& \Theta((\Phi_{\text{gpd},\rtimes}^{-1} \circ \Phi)(f\otimes\1z_{\{(d,c)\}})\bullet \phi)(\zf,(t,a)) \nonumber \\
  &=& ((\Phi_{\text{gpd},\rtimes}^{-1} \circ \Phi)(f\otimes \1z_{\{(d,c)\}}) \bullet \phi)([ta,\zf a],(-ta,a^{-1})) \text{ by } \eqref{def-theta}. \label{left-R-1}
\egln
Moreover, use formula~\eqref{MRW-leftaction} with $\zeta = ([ta,\zf a],(-ta,a^{-1}))$. As we must have $r(\gamma)=r(\zeta)$, the element $\gamma$ in \eqref{MRW-leftaction} is of the form $([ta,\zf a],(s,c'))$. Then
\bgloz
  \gamma^{-1} \zeta = ([(ta+s)c',\zf a c'], -(ta+s)c',(c')^{-1} a^{-1}).
\egloz
Using that $\lambda^{r(\zeta)}=\lambda^{[ta,\zf a]}=\delta_{[ta,\zf a]} \times \lambda_{\Rz \rtimes \spkl{2}}$, we continue the computation at \eqref{left-R-1} and obtain
\bglnoz
  && ((f\otimes \1z_{\{(d,c)\}}) \cdot \Theta(\phi))(\zf,(t,a)) \\
  &=& \int_{\Rz \rtimes \spkl{2}} \lge (\Phi_{\text{gpd},\rtimes}^{-1} \circ \Phi)(f\otimes\1z_{\{(d,c)\}})([ta,\zf a],(s,c')) \cdot \right. \\
  && \left. \cdot \phi([(ta+s)c',\zf a c'], -(ta+s)c',(c')^{-1} a^{-1}) \rge d\lambda(s,c') \\
  &=& c^{\halb} \e(-p(\zf ad)) \e(tad) \int_{\Rz} \e(sd) \check{f}(s) (\Theta(\phi))(\zf,t+sa^{-1},ac)ds.
\eglnoz

The right action of $C_c(\cG_N^N)$ on $C_c(\cG_N)$ is given by
\bgl
\label{MRW-rightaction}
  \phi \bullet g(\zeta) = \int_{\cG_N^N} \phi(\zeta \eta) g(\eta^{-1}) d\lambda^{s(\zeta)}(\eta) \fuer \zeta \in \cG_N.
\egl
It is transformed into the right action $\Theta(\phi) \cdot g = \Theta(\phi \bullet g)$ on $X_0$. By applying $\Psi \circ \Psi_{N, \ti{N}}$, this becomes a right action of $Q = \staralg(u,s_2) \subseteq \cQ_2$ on $X_0$. It suffices to specify the action of each of $u$ and $s_2$. It follows from Corollary~\ref{def-Psi} that
\begin{align}
  s_2 &= (\Psi \circ \Psi_{N, \ti{N}})(\1z_{\pi(\{0\}\times 2\Zz_2) \times \{(0,2^{-1})\}}) \text{ and } \label{s2-from-Psi} \\
  u &= (\Psi \circ \Psi_{N, \ti{N}})(\1z_{N \times \{(1,1)\}}). \label{u-from-Psi}
\end{align}

Take $(\zf,(t,a)) \in \Zz_2 \rtimes (\Rz \rtimes \spkl{2})$. Set $\zeta=([ta,\zf a],(-ta,a^{-1})) \in \cG_N$ in \eqref{MRW-rightaction}. Note that $\eta$ in \eqref{MRW-rightaction} is of the form $([0,\zf],(d,c))$ because of the condition $s(\zeta)=r(\eta)$. Thus we obtain
\bglnoz
(\Theta(\phi) \cdot u)(\zf,(t,a))
  &=& \Theta(\phi \bullet \1z_{N \times \{(1,1)\}})(\zf,(t,a)) \\
  &=& (\phi \bullet \1z_{N \times \{(1,1)\}})([ta,\zf a],(-ta,a^{-1})) \\
  &=& \sum_{(d,c) \in \Zz[\halb] \rtimes \spkl{2}} \lge \phi(([ta,\zf a],(-ta,a^{-1}))([0,\zf],(d,c)))\cdot \right. \\
  && \left. \cdot\1z_{N \times \{(1,1)\}}([-dc,\zf c],(-dc,c^{-1}) \rge \\
  &=& \phi(([ta,\zf a],(-ta,a^{-1}))([0,\zf],(-1,1))) \1z_N([1,\zf]) \\
  &=& \phi(([ta,\zf a],(-ta-a,a^{-1})) \\
  &=& (\Theta(\phi))(\zf+1,t+1,a)
\eglnoz
and
\bglnoz
  (\Theta(\phi) \cdot s_2)(\zf,(t,a))
  &=& (\phi \bullet \1z_{\pi(\gekl{0} \times 2\Zz_2) \times \{(0,2^{-1})\}})([ta,\zf a],(-ta,a^{-1})) \\
  &=& \sum_{(d,c) \in \Zz[\halb] \rtimes \spkl{2}} \lge \phi(([ta,\zf a],(-ta,a^{-1}))([0,\zf],(d,c)))\cdot \right. \\
  && \left. \cdot\1z_{\pi(\gekl{0} \times 2\Zz_2) \times \{(0,2^{-1})\}}([-dc,\zf c],(-dc,c^{-1}) \rge \\
  &=& \phi(([ta,\zf a],(-ta,a^{-1}))([0,\zf],(0,2))) \1z_{\pi(\gekl{0} \times 2\Zz_2)}([0,2\zf]) \\
  &=& \phi(([ta,\zf a],(-ta,a^{-1}2)) \\
  &=& (\Theta(\phi))(2\zf,2t,2^{-1}a).
\eglnoz

Now we turn attention to the inner products. Note that in the proof of Theorem~2.8 in \cite{MRW}, the correct formula for the $C_c(\cG)$-valued inner product should be
\bgl
\label{MRW-leftinnerproduct}
  _{C_c(\cG)}\spkl{\phi_1,\phi_2} (\gamma) = \int_{\cG_N^N} \phi_1(\gamma \zeta \eta) \overline{\phi_2}(\zeta \eta) d\mu^{s(\zeta)}(\eta)
\egl
for $\phi_i\in C_c(\cG_N)$, $i=1,2$, $\gamma \in \cG$ and $\zeta \in \cG_N$ with $r(\zeta)=s(\gamma)$. We then obtain a $C_c(\cG)$-valued
pre-inner product on $X_0$ by
\bgloz
  _{C_c(\cG)}\spkl{\Theta(\phi_1),\Theta(\phi_2)} \defeq \ _{C_c(\cG)}\spkl{\phi_1,\phi_2}
\egloz
for $\phi_i \in C_c(\cG_N)$, $i=1,2$.

Set $\gamma=([r,\xf],(t,a))$ in \eqref{MRW-leftinnerproduct}. We choose
\bgloz
  \zeta = ([(r+t)a,\xf a],(-(r+t)a+p(\xf)a,a^{-1}))
\egloz
so that $r(\zeta)=[(r+t)a,\xf a]$ and $s(\zeta)=[p(\xf),\xf]=[0,\xf-p(\xf)] \in N$. Recall that $p$ was defined in Section~\ref{dualthm}, after \eqref{dual_of_Z2}. Here, we mean by $p(\xf)$ an element of $\Zz[\halb]$ which represents $p(\xf) \in \Zz[\halb] / \Zz$. The condition $r(\eta)=s(\zeta)$ forces $\eta$ to be of the form
\bgl
\label{sumover?}
  ([0,\xf-p(\xf)],(d,c)) \text{ with } s(\eta)=(\xf-p(\xf)-d)c \in \Zz_2.
\egl
With these choices, we obtain
\bgl
\label{gammazetaeta}
  \gamma \zeta \eta = ([r,\xf],(-r+p(\xf)+d,c))
\egl
and
\bgl
\label{zetaeta}
  \zeta \eta = ([(r+t)a,\xf a],(-(r+t)a+p(\xf)a+da,a^{-1}c)).
\egl
Plugging \eqref{gammazetaeta} and \eqref{zetaeta} into \eqref{MRW-leftinnerproduct}, we end up with
\bgln
\label{leftinnerproduct1}
  && _{C_c(\cG)} \spkl{\phi_1,\phi_2}([r,\xf],(t,a)) \\
  &=& \sum_{(d,c)} \lge \phi_1([r,\xf],(-r+p(\xf)+d,c)) \cdot \right. \nonumber \\
  && \left. \cdot \overline{\phi_2}([(r+t)a,\xf a],(-(r+t)a+p(\xf)a+da,a^{-1}c)) \rge \nonumber
\egln
where the sum is taken over all $(d,c)$ in $\Zz[\halb] \rtimes \spkl{2}$ with the property that $(\xf-p(\xf)-d)c$ lies in $\Zz_2$ (see \eqref{sumover?}). Thus we obtain for $\varphi_1$, $\varphi_2$ in $X_0$:
\bgln
\label{R-ip1}
  && _{C_c(\cG)} \spkl{\varphi_1,\varphi_2}([r,\xf],(t,a)) \\
  &=& \sum_{(d,c)} \lge \varphi_1((\xf-p(\xf)-d)c,((r-p(\xf)-d)c,c^{-1})) \cdot \right. \nonumber \\
  && \left. \cdot \overline{\varphi_2}((\xf-p(\xf)-d)c,((r+t)c-p(\xf)c-dc,c^{-1}a)) \rge \nonumber
\egln
by setting $\phi_i=\Theta^{-1}(\varphi_i)$ in \eqref{leftinnerproduct1}. Again the sum is taken over all $(d,c)$ in $\Zz[\halb] \rtimes \spkl{2}$ with the property that $(\xf-p(\xf)-d)c$ lies in $\Zz_2$. Then the required $R$-valued inner product on $X_0$ is obtained as
\bgl
\label{R-ip2}
  _{R} \spkl{\varphi_1,\varphi_2}=(\Phi^{-1} \circ \Phi_{\text{gpd},\rtimes})(_{C_c(\cG)}\spkl{\varphi_1,\varphi_2}).
\egl
We remark that this description of the $R$-valued inner product is not used in the sequel; we have included it for the sake of completeness.

Finally, we obtain a $Q$-valued inner product on $X_0$ by letting
\bgl
\label{Q-ip1}
  \spkl{\varphi_1,\varphi_2}_Q \defeq (\Psi \circ \Psi_{N,\ti{N}})(\spkl{\Theta^{-1}(\varphi_1),\Theta^{-1}(\varphi_2)}_{C_c(\cG_N^N)})
\egl
for $\varphi_1$, $\varphi_2$ in $X_0$. By \eqref{dense-lin-span} it suffices to specify $\spkl{\varphi_1,\varphi_2}_Q$ for
\bglnoz
  && \varphi_1 \defeq \1z_{l_1+k_1\Zz_2} \otimes \xi_1 \otimes \1z_{\gekl{m_1}} \text{ and }\\
  && \varphi_2 \defeq \1z_{l_2+k_2\Zz_2} \otimes \xi_2 \otimes \1z_{\gekl{m_2}},
\eglnoz
with $k_1,k_2 \in [2\rangle$, $l_1,l_2 \in \Zz$, $\xi_1,\xi_2 \in C_c(\Rz)$ and $m_1, m_2 \in \spkl{2}$.

The $C_c(\cG_N^N)$-valued inner product on $C_c(\cG_N)$ is given in the proof of \cite[Theorem~2.8]{MRW} by the formula
\bgl
\label{MRW-rightinnerproduct}
  \spkl{\phi_1,\phi_2}_{C_c(\cG_N^N)}(\eta)=\int_{\cG} \overline{\phi_1}(\gamma^{-1}\zeta)\phi_2(\gamma^{-1}\zeta\eta) d\lambda^{r(\zeta)}(\gamma),
\egl
for $\phi_i \in C_c(\cG_N)$, $i=1,2$, where the element $\zeta$ is chosen in $\cG_N$ so that $s(\zeta)=r(\eta)$.

Set $\eta=([0,\zf],(d,c))$ in \eqref{MRW-rightinnerproduct} and choose $\zeta=([0,\zf],(0,1))$ so that $s(\zeta)=[0,\zf]=r(\eta)$. The condition $r(\gamma)=r(\zeta)=[0,\zf]$ forces $\gamma$ to be of the form $([0,\zf],(t,a))$. With these choices we obtain
\bgl
\label{gamma^-1zeta}
  \gamma^{-1}\zeta=([ta,\zf a],(-ta,a^{-1}))
\egl
and
\bgl
\label{gamma^-1zetaeta}
  \gamma^{-1}\zeta\eta=([ta,\zf a],(-ta+da,a^{-1}c)).
\egl
Plugging  \eqref{gamma^-1zeta} and \eqref{gamma^-1zetaeta} into \eqref{MRW-rightinnerproduct}, we end up with
\bglnoz
  && \spkl{\Theta^{-1}(\varphi_1),\Theta^{-1}(\varphi_2)}_{C_c(\cG_N^N)}([0,\zf],(d,c))) \\
  &=& \int_{\Rz \rtimes \spkl{2}} \lge \overline{\Theta^{-1}(\varphi_1)}([ta,\zf a],(-ta,a^{-1})) \cdot \right. \\
  && \left. \cdot \Theta^{-1} (\phi_2)([ta,\zf a],(-ta+da,a^{-1}c)) \rge d\lambda(t,a) \\
  &=& \int_{\Rz \rtimes \spkl{2}} \overline{\varphi_1}(\zf,(t,a)) \varphi_2((\zf -d)c,(tc-dc,c^{-1}a)) d\lambda(t,a) \\
  &=& \int_{\Rz \rtimes \spkl{2}} \lge \1z_{l_1+k_1\Zz_2}(\zf) \overline{\xi_1}(t) \1z_{\gekl{m_1}}(a) \right. \\
  && \left. \cdot\1z_{l_2+k_2\Zz_2}((\zf -d)c) \xi_2((t-d)c) \1z_{\gekl{m_2}}(c^{-1}a) \rge d\lambda(t,a) \\
  &=& m_1 \delta_{\tfrac{m_1}{m_2},c} \int_\Rz \1z_{l_1+k_1\Zz_2}(\zf) \1z_{l_2+k_2\Zz_2}((\zf -d)\tfrac{m_1}{m_2}) \overline{\xi_1}(t) \xi_2((t-d)\tfrac{m_1}{m_2}) dt \\
  &=& m_1 \delta_{\tfrac{m_1}{m_2},c} \1z_{l_1+k_1\Zz_2}(\zf) \1z_{l_2+k_2\Zz_2}((\zf -d)\tfrac{m_1}{m_2}) \spkl{\xi_1,\xi_2((\sqcup-d)\tfrac{m_1}{m_2})}_{\Ell^2(\Rz)}
\eglnoz
where $\spkl{\xi_1,\xi_2}_{\Ell^2(\Rz)} = \int_{\Rz} \overline{\xi_1}(t) \xi_2(t) dt$.

We distinguish two cases:

{\bf Case (1): $c \leq 1$, i.e. $m_1 \leq m_2$}. We claim that
\bgln
\label{Q-ip2}
  && (\Psi \circ \Psi_{N,\ti{N}})(\spkl{\Theta^{-1}(\varphi_1),\Theta^{-1}(\varphi_2)}_{C_c(\cG_N^N)}) \\
  &=& m_1 \sum_{b \in \Zz} \spkl{\xi_1,\xi_2((\sqcup+b)\tfrac{m_1}{m_2})}_{\Ell^2(\Rz)} u^{l_1} e_{k_1} u^{-l_1} u^{-b} s_{\tfrac{m_2}{m_1}} u^{l_2} e_{k_2} u^{-l_2}.
  \nonumber
\egln
To show this, set $\Lambda_b = m_1 \spkl{\xi_1,\xi_2((\sqcup+b)\tfrac{m_1}{m_2})}_{\Ell^2(\Rz)}$ for every $b \in \Zz$. Using
that  $u=\Psi(\1z_{\Zz_2 \times \gekl{(1,1)}})$ by \eqref{u-from-Psi} and $s_2=\Psi(\1z_{2\Zz_2 \times \gekl{(0,2^{-1})}})$ by \eqref{s2-from-Psi},
we obtain
\bglnoz
  && u^{l_1} e_{k_1} u^{-l_1}=u^{l_1} s_{k_1} s_{k_1}^* u^{-l_1} \\
  &=& \Psi
  \rukl{\1z_{\Zz_2 \times \gekl{(l_1,1)}}*\1z_{k_1\Zz_2 \times \gekl{(0,k_1^{-1})}}*(\1z_{k_1\Zz_2 \times \gekl{(0,k_1^{-1})}})^**\1z_{\Zz_2 \times \gekl{(-l_1,1)}}}
\eglnoz
where $*$ denotes the product in $C^*(\ti{\cG}_{\ti{N}}^{\ti{N}})$.

A routine (but long) computation in $C^*(\ti{\cG}_{\ti{N}}^{\ti{N}})$, using the choice of Haar system on $\ti{\cG}_{\ti{N}}^{\ti{N}}$ given by counting measure on  $\Zz[\halb] \rtimes \spkl{2}$, shows that
\bgloz
  \1z_{(l+k\Zz_2)\times\gekl{(0,1)}}=\1z_{\Zz_2 \times \gekl{(l,1)}}*\1z_{k\Zz_2 \times \gekl{(0,k^{-1})}}*(\1z_{k\Zz_2 \times \gekl{(0,k^{-1})}})^**\1z_{\Zz_2 \times \gekl{(-l,1)}}
\egloz
for all $l\in \Zz$ and $k\in \spkl{2}$. We leave the details to the reader.

A few more computations show that
\bglnoz
  && \1z_{\Zz_2 \times \gekl{(-b,1)}}*\1z_{\tfrac{m_2}{m_1}\Zz_2 \times \gekl{(0,\tfrac{m_1}{m_2})}}*\1z_{(l_2+k_2\Zz_2)\times\gekl{(0,1)}}(\zf,(d,c)) \\
  &=& \1z_{\Zz_2}(\zf) \1z_{l_2+k_2\Zz_2}((\zf+b)\tfrac{m_1}{m_2}) \1z_{\gekl{-b}}(d) \1z_{\gekl{\tfrac{m_1}{m_2}}}(c)
\eglnoz
and thus that
\bgloz
  \1z_{(l_1+k_1\Zz_2)\times\gekl{(0,1)}}*\1z_{\Zz_2 \times \gekl{(-b,1)}}*\1z_{\tfrac{m_1}{m_2}\Zz_2 \times \gekl{(0,\tfrac{m_1}{m_2})}}
  *\1z_{(l_2+k_2\Zz_2)\times\gekl{(0,1)}}
\egloz
evaluated at $(\zf,(d,c))$ gives
\bgloz
  \1z_{l_1+k_1\Zz_2}(\zf) \1z_{l_2+k_2\Zz_2}((\zf+b)\tfrac{m_1}{m_2}) \1z_{\gekl{-b}}(d) \1z_{\gekl{\tfrac{m_1}{m_2}}}(c).
\egloz
Putting all this together, we get
\bglnoz
  && (\Psi \circ \Psi_{N,\ti{N}})^{-1}(\sum_b \Lambda_b u^{l_1} e_{k_1} u^{-l_1} u^{-b} s_{\tfrac{m_1}{m_2}} u^{l_2} e_{k_2} u^{-l_2})(\zf,(d,c)) \\
  &=& \sum_b \Lambda_b \lge \1z_{(l_1+k_1\Zz_2)\times\gekl{(0,1)}}*\1z_{\Zz_2 \times \gekl{(-b,1)}}*\1z_{\tfrac{m_1}{m_2}\Zz_2 \times \gekl{(0,\tfrac{m_1}{m_2})}} \right.
  \\
  && \left. *\1z_{(l_2+k_2\Zz_2)\times\gekl{(0,1)}} \rge(\zf,(d,c)) \\
  &=& \sum_b \Lambda_b \1z_{l_1+k_1\Zz_2}(\zf) \1z_{l_2+k_2\Zz_2}((\zf+b)\tfrac{m_1}{m_2}) \1z_{\gekl{-b}}(d) \1z_{\gekl{\tfrac{m_1}{m_2}}}(c) \\
  &=& \Lambda_{-d} \delta_{\tfrac{m_1}{m_2},c} \1z_{l_1+k_1\Zz_2}(\zf) \1z_{l_2+k_2\Zz_2}((\zf-d)\tfrac{m_1}{m_2}) \\
  &=& \spkl{\Theta^{-1}(\varphi_1),\Theta^{-1}(\varphi_2)}_{C_c(\cG_N^N)}.
\eglnoz
This proves our claim and thus settles the first case.

{\bf Case (2): $c \geq 1$, i.e. $m_1 \geq m_2$}. Proceeding analogously as in the first case, we can compute that
\bgln
\label{Q-ip3}
  && (\Psi \circ \Psi_{N,\ti{N}})(\spkl{\Theta^{-1}(\varphi_1),\Theta^{-1}(\varphi_2)}_{C_c(\cG_N^N)}) \\
  &=& m_1 \sum_{b \in \Zz} \spkl{\xi_1,\xi_2((\sqcup+b)\tfrac{m_1}{m_2})}_{\Ell^2(\Rz)}
  u^{l_1} e_{k_1} u^{-l_1} s_{\tfrac{m_1}{m_2}}^* u^{-\tfrac{m_1}{m_2}b} u^{l_2} e_{k_2} u^{-l_2}. \nonumber
\egln

We have thus proven
\bprop
\label{X_0}
The space $X_0 = C_c(\Zz_2 \times (\Rz \rtimes \spkl{2}))$ carries the following structure of an $R$--$Q$ pre-imprimitivity bimodule:

The left $R$-action is given by
\bgln
  \label{leftaction}
  && ((f \otimes \1z_{\gekl{(d,c)}}) \cdot \varphi)(\zf,(t,a)) \\
  &=& c^{\halb} \e(-p(\zf ad)) \e(tad) \int_{\Rz} \e(sd) \check{f}(s) \varphi(\zf,t+sa^{-1},ac)ds \nonumber
\egln
and the right $Q$-action is given by
\begin{align}
  (\varphi\cdot u)(\zf,(t,a))&=\varphi(\zf+1,t+1,a),\label{u-from-right}\\
  (\varphi\cdot s_2)(\zf,(t,a))&=\varphi(2\zf,2t,2^{-1}a)\label{s2-from-right}
\end{align}
for all $f \otimes \1z_{\gekl{(d,c)}} \in R$, $\varphi \in X_0$ and $(\zf,(t,a)) \in \Zz_2 \times (\Rz \rtimes \spkl{2})$;
the left $R$-valued inner product is given by \eqref{R-ip2} together with \eqref{R-ip1},
and the right $Q$-valued inner product is given by \eqref{Q-ip1} together with \eqref{Q-ip2} and \eqref{Q-ip3}.
\eprop

\section{Induced representations}
Recall that $\lambda_2$ is the (restricted) left regular representation of $\cQ_2$ on $\ell^2(\Zz)$,
and that $\pi$ is the standard representation of $C_0(\Rz) \rtimes \Zz[\halb] \rtimes \spkl{2}$ on $\Ell^2(\Rz)$ from \eqref{standard-rep}.
Let $X$ denote the completion of $X_0$ from Proposition~\ref{X_0} to a $(C_0(\Rz) \rtimes \Zz[\halb] \rtimes \spkl{2})$--$\cQ_2$-imprimitivity bimodule. The goal of this section is to establish the following theorem:

\btheo
\label{lambda2-equiv-pi}
The representations $\pi$ and $X \text{--} \Ind (\lambda_2)$ are unitarily equivalent.
\etheo

This is our third main result. Its proof consists of two parts. First, we construct a unitary isomorphism of the Hilbert spaces $\Ell^2(\Rz)$ and $X \otimes_{\cQ_2, \lambda_2} \ell^2(\Zz)$. The second Hilbert space is precisely the one on which $X$--$\Ind (\lambda_2)$ is defined. Let
\bgloz
  j: X_0 \odot_{Q, \lambda_2} \ell^2(\Zz) {\ri} X \otimes_{\cQ_2, \lambda_2} \ell^2(\Zz)
\egloz
be given by $\varphi \otimes \varepsilon \ma \varphi^{\bullet} \otimes \varepsilon$ where $X_0 \ni \varphi \ma \varphi^{\bullet} \in X$ is the canonical map.

\blemma
The map
\bgloz
  C_c(\Rz) \ni \xi \overset{W}{\ma} j((\1z_{\Zz_2} \otimes \xi \otimes \1z_{\gekl{1}}) \otimes \varepsilon_0) \in X \otimes_{\cQ_2, \lambda_2} \ell^2(\Zz)
\egloz
extends to a unitary isomorphism $\Ell^2(\Rz) \overset{\cong}{\ri} X \otimes_{\cQ_2, \lambda_2} \ell^2(\Zz)$.
\elemma
Here $\menge{\varepsilon_n}{n \in \Zz}$ is the canonical orthonormal basis of $\ell^2(\Zz)$, so $\varepsilon_0$ is the $0$-th basis vector.
\bproof
We prove first that $W$ preserves the inner products. Take $\xi_1$, $\xi_2$ in $C_c(\Rz)$. We have, using \eqref{Q-ip2},
\bglnoz
  && \spkl{W(\xi_1),W(\xi_2)}_{X \otimes_{\cQ_2, \lambda_2} \ell^2(\Zz)} \\
  &=& \spkl{j((\1z_{\Zz_2} \otimes \xi_1 \otimes \1z_{\gekl{1}}) \otimes \varepsilon_0),
  j((\1z_{\Zz_2} \otimes \xi_2 \otimes \1z_{\gekl{1}}) \otimes \varepsilon_0)}_{X \otimes_{\cQ_2, \lambda_2} \ell^2(\Zz)} \\
  &=& \spkl{\varepsilon_0,
  \lambda_2(\spkl{\1z_{\Zz_2} \otimes \xi_1 \otimes \1z_{\gekl{1}},\1z_{\Zz_2} \otimes \xi_2 \otimes \1z_{\gekl{1}}}_{\cQ_2}) \varepsilon_0}_{\ell^2(\Zz)} \\
  &=& \spkl{\varepsilon_0,\sum_{b \in \Zz} \spkl{\xi_1,\xi_2(\sqcup+b)}_{\Ell^2(\Rz)} \lambda_2(u^{-b}) \varepsilon_0}_{\ell^2(\Zz)} \\
  &=& \spkl{\varepsilon_0,\sum_{b \in \Zz} \spkl{\xi_1,\xi_2(\sqcup+b)}_{\Ell^2(\Rz)} \varepsilon_{-b}}_{\ell^2(\Zz)} \\
  &=& \spkl{\varepsilon_0,\spkl{\xi_1,\xi_2}_{\Ell^2(\Rz)} \varepsilon_0}_{\ell^2(\Zz)} = \spkl{\xi_1,\xi_2}_{\Ell^2(\Rz)}.
\eglnoz
Thus $W$ extends to an isometry $\Ell^2(\Rz) \ri X \otimes_{\cQ_2, \lambda_2} \ell^2(\Zz)$, and it remains to prove that $\img(W)$ is dense in $X \otimes_{\cQ_2, \lambda_2} \ell^2(\Zz)$. By \eqref{dense-lin-span}
\bgloz
  \menge{j((\1z_{l+k\Zz_2} \otimes \xi \otimes \1z_{\gekl{m}}) \otimes \varepsilon_n)}{k \in [2\rangle, l \in \Zz, \xi \in C_c(\Rz), m \in \spkl{2}, n \in \Zz}
\egloz
is total in $X \otimes_{\cQ_2, \lambda_2} \ell^2(\Zz)$. As $\img(W)$ is a subspace of $X \otimes_{\cQ_2, \lambda_2} \ell^2(\Zz)$, it suffices to prove that for all $k \in [2\rangle$, $l \in \Zz$, $\xi \in C_c(\Rz)$, $m \in \spkl{2}$ and $n \in \Zz$, $j((\1z_{l+k\Zz_2} \otimes \xi \otimes \1z_{\gekl{m}}) \otimes \varepsilon_n)$ lies in $\img(W)$.

We have
\bglnoz
  && j((\1z_{l+k\Zz_2} \otimes \xi \otimes \1z_{\gekl{m}}) \otimes \varepsilon_n) = (\1z_{l+k\Zz_2} \otimes \xi \otimes \1z_{\gekl{m}})^{\bullet} \otimes \varepsilon_n \\
  &=& (\1z_{l+k\Zz_2} \otimes \xi \otimes \1z_{\gekl{m}})^{\bullet} \otimes \lambda_2(u^n) \varepsilon_0 \\
  &=& ((\1z_{l+k\Zz_2} \otimes \xi \otimes \1z_{\gekl{m}}) \cdot u^n)^{\bullet} \otimes \varepsilon_0 \\
  &=& ((\1z_{-n+l+k\Zz_2} \otimes \xi(\sqcup+n) \otimes \1z_{\gekl{m}}))^{\bullet} \otimes \varepsilon_0 \text{ by }\eqref{s2-from-right}.
\eglnoz
Thus it suffices to prove that for all $k \in [2\rangle$, $l \in \Zz$, $\xi \in C_c(\Rz)$ and $m \in \spkl{2}$, $j((\1z_{l+k\Zz_2} \otimes \xi \otimes \1z_{\gekl{m}}) \otimes \varepsilon_0)$ lies in $\img(W)$.

As a next step, we claim that
\bgl
\label{l+kZvanishes}
  j((\1z_{l+k\Zz_2} \otimes \xi \otimes \1z_{\gekl{m}}) \otimes \varepsilon_0) = 0 \text{ in } X \otimes_{\cQ_2, \lambda_2} \ell^2(\Zz) \text{ unless } l \in k\Zz_2.
\egl
To see this, use \eqref{Q-ip2} and compute the inner product of $j((\1z_{l+k\Zz_2} \otimes \xi \otimes \1z_{\gekl{m}}) \otimes \varepsilon_0)$ with itself in $X \otimes_{\cQ_2, \lambda_2} \ell^2(\Zz)$. We obtain
\bgloz
  m \spkl{\varepsilon_0,\sum_{b \in \Zz} \spkl{\xi,\xi(\sqcup+b)}_{\Ell^2(\Rz)} \lambda_2(u^l e_k u^{-l} u^{-b} u^l e_k u^{-l}) \varepsilon_0}_{\ell^2(\Zz)}.
\egloz
Since $\lambda_2(u^l e_k u^{-l}) \varepsilon_n=\1z_{l+k\Zz}(n)\varepsilon_n$, we see that $\lambda_2(u^l e_k u^{-l}) \varepsilon_0 \neq 0$ only when $0 \in l+k\Zz$. The last condition is equivalent to $l\in k\Zz$, and to $l \in k\Zz_2$. Therefore we only have to prove that for all $k \in [2\rangle$, $\xi \in C_c(\Rz)$ and $m \in \spkl{2}$, $j((\1z_{k\Zz_2} \otimes \xi \otimes \1z_{\gekl{m}}) \otimes \varepsilon_0)$ lies in $\img(W)$.

Once again using \eqref{Q-ip2}, we obtain that the inner product of
\bgloz
  j((\1z_{k \Zz_2} \otimes \xi \otimes \1z_{\gekl{m}}) \otimes \varepsilon_0) - j((\1z_{\Zz_2} \otimes \xi \otimes \1z_{\gekl{m}}) \otimes \varepsilon_0)
\egloz
with itself in $X \otimes_{\cQ_2, \lambda_2} \ell^2(\Zz)$ can be written as
\bglnoz
  && m \spkl{\varepsilon_0,\spkl{\xi,\xi}_{\Ell^2(\Rz)} \lambda_2(e_k^2) \varepsilon_0}_{\ell^2(\Zz)}
  - m \spkl{\varepsilon_0,\spkl{\xi,\xi}_{\Ell^2(\Rz)} \lambda_2(e_k) \varepsilon_0}_{\ell^2(\Zz)} \\
  && - m \spkl{\varepsilon_0,\spkl{\xi,\xi}_{\Ell^2(\Rz)} \lambda_2(e_k) \varepsilon_0}_{\ell^2(\Zz)}
  + m \spkl{\varepsilon_0,\spkl{\xi,\xi}_{\Ell^2(\Rz)} \varepsilon_0}_{\ell^2(\Zz)} = 0.
\eglnoz
This implies
\bgl
\label{indHsp1}
  j((\1z_{k\Zz_2} \otimes \xi \otimes \1z_{\gekl{m}}) \otimes \varepsilon_0) = j((\1z_{\Zz_2} \otimes \xi \otimes \1z_{\gekl{m}}) \otimes \varepsilon_0).
\egl
Thus it remains to show that for all $\xi \in C_c(\Rz)$ and $m \in \Zz$, the element $j((\1z_{\Zz_2} \otimes \xi \otimes \1z_{\gekl{m}}) \otimes \varepsilon_0)$ lies in $\img(W)$.

Assume first that $m \in [2\rangle$, i.e. $\log_2(m) \geq 0$. Then \eqref{s2-from-right} implies
\bgln
  \label{indHsp2}
  && j((\1z_{\Zz_2} \otimes \xi \otimes \1z_{\gekl{m}}) \otimes \varepsilon_0) = j((\1z_{\Zz_2} \otimes \xi(\sqcup \cdot m^{-1}) \otimes \1z_{\gekl{1}}) \cdot s_m \otimes \varepsilon_0) \\
  &=& j((\1z_{\Zz_2} \otimes \xi( \sqcup \cdot m^{-1}) \otimes \1z_{\gekl{1}}) \otimes \lambda_2(s_m) \varepsilon_0) \nonumber \\
  &=& j((\1z_{\Zz_2} \otimes \xi( \sqcup \cdot m^{-1}) \otimes \1z_{\gekl{1}}) \otimes \varepsilon_0)
  = W(\xi( \sqcup \cdot m^{-1})), \nonumber
\egln
which lies in $\img(W)$.

If $m^{-1}$ lies in $[2\rangle$, i.e. $\log_2(m) \leq 0$, then by \eqref{s2-from-right}
\bglnoz
  && j((\1z_{\Zz_2} \otimes \xi \otimes \1z_{\gekl{m}}) \otimes \varepsilon_0) = j((\1z_{\Zz_2} \otimes \xi \otimes \1z_{\gekl{m}}) \otimes
  \lambda_2(s_{m^{-1}}) \varepsilon_0) \\
  &=& j((\1z_{\Zz_2} \otimes \xi \otimes \1z_{\gekl{m}}) \cdot s_{m^{-1}} \otimes \varepsilon_0) \\
  &=& j((\1z_{\Zz_2} \otimes \xi( \sqcup \cdot m^{-1}) \otimes \1z_{\gekl{1}}) \otimes \varepsilon_0)
  = W(\xi( \sqcup \cdot m^{-1})),
\eglnoz
which lies in $\img(W)$.

This proves that $\img(W)$ is dense in $X \otimes_{\cQ_2, \lambda_2} \ell^2(\Zz)$. Thus $W$ extends to a unitary isomorphism $\Ell^2(\Rz) \overset{\cong}{\ri} X \otimes_{\cQ_2, \lambda_2} \ell^2(\Zz)$, as claimed.
\eproof

We denote the extension of $W$ to $\Ell^2(\Rz)$ by $W$ as well, and we let $F$ denote the Fourier transform, seen as a unitary on $\Ell^2(\Rz)$. Thus, by our convention, $(F\xi)(t)=\int_\Rz \e(tx)\xi(x)dx$ for $\xi\in C_c(\Rz)\subseteq \Ell^2(\Rz)$.

The final step in proving Theorem~\ref{lambda2-equiv-pi} is to prove the following result:
\blemma
We have $
  \Ad(W^{-1}) \circ (X \text{--} \Ind(\lambda_2)) = \Ad(F) \circ \pi$.
\elemma
\bproof
For an element $f \otimes \1z_{\gekl{(d,c)}} \in R \subseteq C_0(\Rz) \rtimes \Zz[\halb] \rtimes \spkl{2}$ with $f$ in $F_{\Rz}(C_c(\Rz))$ and for $\xi \in C_c(\Rz) \subseteq \Ell^2(\Rz)$ denote by $\eta_{\xi}^{f,d,c}$ the function on $\Rz$ given by
$$
  \eta_{\xi}^{f,d,c}(t) \defeq \e(tc^{-1}d) \int_{\Rz} \e(sd) \check{f}(s) \xi(t+sc) ds.
$$
Inserting $\varphi = \1z_{\Zz_2} \otimes \xi \otimes \1z_{\gekl{1}}$ into \eqref{leftaction}, we obtain
\bglnoz
  && (f \otimes \1z_{\gekl{(d,c)}}) \cdot (\1z_{\Zz_2} \otimes \xi \otimes \1z_{\gekl{1}})(\zf,(t,a)) \\
  &=& c^{\halb} \e(-p(\zf c^{-1}d)) \e(tc^{-1}d) \int_{\Rz} \e(sd) \check{f}(s) \xi(t+sc) ds \ \1z_{\gekl{c^{-1}}}(a),
\eglnoz
thus
\bgloz
  (f \otimes \1z_{\gekl{(d,c)}}) \cdot (\1z_{\Zz_2} \otimes \xi \otimes \1z_{\gekl{1}})
  = c^{\halb} \e(-p(\sqcup c^{-1}d)) \otimes \eta_\xi^{f,d,c} \otimes \1z_{\gekl{c^{-1}}}.
\egloz
Then we compute
\bgln
  && W^{-1} (X \text{--} \Ind(\lambda_2)) (f \otimes \1z_{\gekl{(d,c)}}) W(\xi) \nonumber \\
  &=& W^{-1} (X \text{--} \Ind(\lambda_2)) (f \otimes \1z_{\gekl{(d,c)}}) (j(\1z_{\Zz_2} \otimes \xi \otimes \1z_{\gekl{1}}) \otimes \varepsilon_0) \nonumber \\
  &=& W^{-1} j((c^{\halb} \e(-p(\sqcup c^{-1}d)) \otimes \eta^{f,d,c}_{\xi} \otimes \1z_{\{c^{-1}\}}) \otimes \varepsilon_0) \label{W*W}.
\egln
To compute $j((c^{\halb} \e(-p(\sqcup c^{-1}d)) \otimes \eta^{f,d,c}_{\xi} \otimes \1z_{\{d^{-1}\}}) \otimes \varepsilon_0)$, we first write the function $\e(-p(\sqcup c^{-1}d))$ as a linear combination of functions of the form $\1z_{l+k\Zz_2}$. Then, we see by \eqref{l+kZvanishes} that only the summand with $l \in k\Zz_2$ gives a contribution.

Writing $c^{-1}d = nm^{-1}$ with $m \in [2\rangle$ and an odd integer $n$ (i.e. $m$ and $n$ are relatively prime), we see that $\e(-p(\zf c^{-1}d))=1$ if $\zf$ lies in $m\Zz_2$. So the function $\e(-p(\sqcup c^{-1}d))$ is constant on cosets of the form $l+m\Zz_2$, and we obtain
\bgloz
  \e(-p(\sqcup c^{-1}d))=\sum_{l=0}^{m-1} \e(-lnm^{-1}) \1z_{l+m\Zz_2},
\egloz
and \eqref{W*W} becomes
\bglnoz
  && W^{-1} (X \text{--} \Ind(\lambda_2)) (f \otimes \1z_{\gekl{(d,c)}}) W(\xi) \\
  &=& W^{-1} j(c^{\halb}((\sum_{l=0}^{m-1} \e(-lnm^{-1}) \1z_{l+m\Zz_2}) \otimes \eta_{\xi}^{f,d,c} \otimes \1z_{\gekl{c^{-1}}}) \otimes \varepsilon_0) \\
  &=& W^{-1} j(\1z_{\Zz_2} \otimes c^{\halb} \eta_{\xi}^{f,d,c} \otimes \1z_{\gekl{c^{-1}}}) \otimes \varepsilon_0) \text{ by }\eqref{indHsp1} \\
  &=& W^{-1} j((\1z_{\Zz_2} \otimes (D_c^* \eta^{f,d,c}_{\xi}) \otimes \1z_{\gekl{1}}) \otimes \varepsilon_0) \text{ by } \eqref{indHsp2} \\
  &=& W^{-1} W(D_c^*\eta^{f,d,c}_{\xi}).
\eglnoz

We claim that $D_c^*\eta^{f,d,c}_{\xi}=F \pi(f \otimes \1z_{\gekl{(d,c)}}) F^{-1} (\xi)$. Indeed,
\bglnoz
  D_c^*\eta^{f,d,c}_{\xi}(t)&=&
  c^{\halb} \e(td) \int_{\Rz} \e(sd) \check{f}(s) \xi((t+s)c) ds \\
  &=& c^{\halb} \int_{\Rz} \e((t+s)d) \check{f}(s) \xi((t+s)c) ds \\
  &=& c^{\halb} \int_{\Rz} \e(sd) \check{f}(-t+s) \xi(sc) ds \\
  &=& c^{\halb} c^{-1} \int_{\Rz} \e(sc^{-1}d) \check{f}(-t+sc^{-1}) \xi(s) ds \\
  &=& c^{-\halb} \int_{\Rz} \int_{\Rz} \e(sc^{-1}d) \e(tr) \e(-sc^{-1}r) f(r) \xi(s) dr ds \\
  &=& c^{-\halb} \int_{\Rz} \e(tr) f(r) (\int_{\Rz} \e(-c^{-1}(r-d)s) \xi(s) ds) dr \\
  &=& (F M_f T_d D_c F^{-1})(\xi)(t) = (F \pi(f \otimes \1z_{\gekl{(d,c)}}) F^{-1})(\xi)(t).
\eglnoz
Since the elements $f \otimes \1z_{\gekl{(d,c)}} \in R$ form a total set in $C_0(\Rz) \rtimes \Zz[\halb] \rtimes \spkl{2}$, we have shown that
\bgloz
  \Ad(W^{-1}) \circ (X \text{--} \Ind(\lambda_2)) = \Ad(F) \circ \pi.
\egloz
\eproof
With these two lemmas, we have proven that the canonical representations correspond to one another with respect to the imprimitivity bimodule from Proposition~\ref{X_0}. In other words, we have proven Theorem~\ref{lambda2-equiv-pi}.

Roughly speaking, this last theorem says that although our proof of the 2-adic duality theorem does not yield an explicit isomorphism, we do obtain an explicit imprimitivity bimodule which arises naturally and which preserves the canonical representations of our algebras.

Let us close with the remark that Theorem~\ref{lambda2-equiv-pi} has its canonical analogue with $\cQ_\Nz$ (or $\cQ_\Zz$) in place of $\cQ_2$. This means that the proof of the duality result \cite[Corollary~3.10]{CuLi2} for $\Gamma=\Qz \pos$ (or $\Qz \reg$) provides a concrete imprimitivity bimodule which implements the Morita equivalence $C_0(\Rz) \rtimes \Qz \rtimes \Qz \pos \sim_M \cQ_\Nz$ (or $C_0(\Rz) \rtimes \Qz \rtimes \Qz \reg \sim_M \cQ_\Zz$) and that the canonical representations are transported into one another by this particular imprimitivity bimodule. The proof of this result for $\cQ_\Nz$ (or $\cQ_\Zz$) should be analogous to the one of Theorem~\ref{lambda2-equiv-pi}.


\begin{thebibliography}{99}

\bibitem[BraJor]{BraJor} O. \textsc{Bratteli} and P. E. T. \textsc{Jorgensen},
	\emph{Isometries, shifts, Cuntz algebras and multiresolution analyses of scale N},
	Int. Eq. \& Op. Th. {\bf 28} (1997), 382–-443.

%\bibitem[BanHLR]{BanHLR} N. \textsc{Brownlowe}, A. \textsc{an Huef}, M. \textsc{Laca} and I.  \textsc{Raeburn},
%  \emph{Boundary quotients of the Toeplitz algebra of the affine semigroup over the natural numbers}, Ergodic Theory \& Dynam. Systems, to appear.

\bibitem[Cun1]{Cun1}J. \textsc{Cuntz}, \emph{Simple $C^*$-algebras generated by
isometries}, Comm. Math. Phys. {\bf 57} (1977), 173--185.

\bibitem[Cun2]{Cun2} J. \textsc{Cuntz}, $C^*$-algebras associated with the $ax + b$-semigroup over $\Nz$, \emph{K-theory and noncommutative geometry}, 201--215, EMS Ser. Congr. Rep., Eur. Math. Soc., Z\"{u}rich, 2008.

\bibitem[CuLi1]{CuLi1} J. \textsc{Cuntz} and X. \textsc{Li},
  {The Regular $C^*$-algebra of an Integral Domain}, \emph{Quanta of Maths}, 149--170, Clay Math. Proc., 11, Amer. Math. Soc., Providence, RI, 2010.

\bibitem[CuLi2]{CuLi2} J. \textsc{Cuntz} and X. \textsc{Li},
  \emph{$C^*$-algebras associated with integral domains and crossed products by actions on adele spaces}, J. Noncommut. Geom. {\bf 5} (2011), 1--37.

\bibitem[Dut]{Dut} D. E. \textsc{Dutkay}, \emph{Low-pass filters and representations of the Baumslag Solitar
group}, Trans. Amer. Math. Soc. {\bf 358} (2006), 5271--5291.

\bibitem[EanHR]{EanHR} R. \textsc{Exel}, A. \textsc{an Huef}, I. \textsc{Raeburn}
  \emph{Purely infinite $C^*$-algebras associated to integer dilation matrices}, Indiana Univ. Math. J., to appear.

\bibitem[Exel]{Exel} R. \textsc{Exel}, \emph{A new look at the crossed product of a $C^*$-algebra by an endomorphism},
Ergodic  Theory \& Dynan. Systems {\bf 23} (2003), 1733--1750.

\bibitem[HeRo]{HeRo} E. \textsc{Hewitt} and K. A. \textsc{Ross}, Abstract Harmonic Analysis, Vol. I and II,
Springer, Berlin, 1963 and 1970.

\bibitem[Hir]{Hir} I. \textsc{Hirshberg},
  \emph{On $C^*$-algebras associated to certain endomorphisms of discrete groups},
  New York J. Math. {\bf 8} (2002), 99--109.

%\bibitem[HLS]{HLS} J. H. \textsc{Hong}, N. S.  \textsc{Larsen} and W.  \textsc{Szymanski},
 % \emph{The Cuntz algebra $\cQ_\Nz$ and $C^*$-algebras of product systems},
%preprint, arXiv:1108.0790[math.OA].

\bibitem[Ka]{Ka} T. \textsc{Katsura}, \emph{A  class of $C^*$-algebras generalizing both graph
algebras and homeomorphism $C^*$-algebras IV, pure infiniteness},  J. Funct. Anal. {\bf 254} (2008), 1161--1187.

\bibitem[Laca]{Laca} M. \textsc{Laca},
  \emph{From endomorphisms to automorphisms and back: dilations and full corners},
  J. London Math. Soc. {\bf 61} (2000), 893-904.

\bibitem[LacRae]{LacRae} M. \textsc{Laca} and I. \textsc{Raeburn},
  \emph{Semigroup crossed products and the Toeplitz algebras of nonabelian groups},
  J. Funct. Anal. {\bf 139} (1996), 415--440.

\bibitem[LarLi]{LarLi}  N. S. \textsc{Larsen} and X. \textsc{Li},
  \emph{Dilations of semigroup crossed products
as crossed products of dilations}, Proc. Amer. Math. Soc., to appear.

\bibitem[LarRae1]{LarRae1} N. S. \textsc{Larsen} and I. \textsc{Raeburn},
	\emph{From filters to wavelets via direct limits},
	in Operator Theory, Operator Algebras and Applications, Contemp. Math., vol. 414, Amer. Math. Soc., Providence, 2006, 35–-40.

\bibitem[LarRae2]{LarRae2} N. S. \textsc{Larsen} and I. \textsc{Raeburn},
 \emph{Projective multi-resolution analyses arising from direct limits of Hilbert modules},
 Math. Scand. {\bf 100} (2007), 317--361.

\bibitem[Li]{Li} X. \textsc{Li},
	\emph{Ring $C^*$-algebras},
	Math. Ann. {\bf 348} (2010), 859--898.

\bibitem[MRW]{MRW} P. S. \textsc{Muhly}, J. N. \textsc{Renault} and D. P. \textsc{Williams},
	\emph{Equivalence and isomorphism for groupoid $C^*$-algebras},
	J. Operator Theory {\bf 17} (1987), 3-22.
	
\bibitem[Pims]{Pims} M. V. \textsc{Pimsner},
  \emph{A class of $C^*$-algebras generalizing both {C}untz-{K}rieger algebras and crossed products by ${\bf Z}$},
  Free probability theory (Waterloo, ON, 1995), Amer. Math. Soc., Providence, RI, 1997, 189--212.

\bibitem[R\o r]{Ror} M. \textsc{R\o rdam},
	\emph{Classification of Nuclear $C^*$-Algebras} in Classification of Nuclear $C^*$-Algebras. Entropy in Operator Algebras,
	Encyclopaedia of Mathematical Sciences, Vol. 126, Springer-Verlag, Berlin Heidelberg New York, 2002.

\bibitem[Wil]{Wil} D. P. \textsc{Williams},
  \emph{Crossed products of $C^*$-algebras},
  Mathematical Surveys and Monographs, vol. 134, Amer. Math. Soc., Providence, RI, 2007.

\end{thebibliography}
\end{document}